\input amstex 
\documentstyle{amsppt} 
\nologo 
\magnification=1200 
\pageheight{18.3cm} 
 
\loadbold 
 
\define\inn#1{{\overset{\lower1.5pt\hbox{$\ssize\circ$}}\to#1}} 
 
\define\eps{\varepsilon} 
\redefine\phi{\varphi} 
\define\k{{\frak k}}

\define\LL{{\Cal L}} 
 
\define\RR{{\Cal R}}

\define\g{{\frak g}} 
\define\h{{\frak h}} 
 
\define\z{{\frak z}}

\define\uu{{\frak u}}

\define\HH{{\Cal H}} 
 
\define\la{\lambda} 
\define\lap{\lambda'} 
\define\nup{\nu'} 
\define\trans{{}^T}

\define\subla{_\lambda} 
 
\define\sublap{_{\lambda'}} 
 
\define\om{\omega} 
\define\Om{\Omega}

\define\scp{{\<\,\,,\,\>}} 
\define\inv{^{-1}} 
\define\Sym{\text{\rm{Sym}}}

\define\Aut{{\text{\rm{Aut}}}} 
\define\Autm{\Aut{}_{g_0}^{\,T}(M)} 
\define\Autphm{\Aut{}_{\phi g_0}^{\,T}(M)} 
\define\Autbm{\overline{\Aut}{}_{g_0}^{\,\,T}(M)} 
\define\Autphbm{\overline{\Aut}{}_{\phi g_0}^{\,\,T}(M)} 
 
\define\Ad{{\text{\rm{Ad}}}} 
\define\ad{{\text{\rm{ad}}}} 
\define\Id{\text{Id}} 
\define\U{{\text{\rm{U}}}} 
\define\SU{{\text{\rm{SU}}}} 
\define\su{{\text{$\frak s$}\text{$\frak u$}}} 
\redefine\O{{\text{\rm{O}}}} 
\define\SO{{\text{\rm{SO}}}} 
\define\so{{\text{$\frak s$}\text{$\frak o$}}} 
\define\SL{{\text{\rm{SL}}}} 
\define\Spin{{\text{\rm{Spin}}}} 
 
\define\dvol{{\text{{\it dvol}}}}

\define\kernn{{\text{\rm{ker}}\,}}

\redefine\exp{{\text{\rm{exp}}}}

\define\Isom{{\text{\rm{Isom}}}} 
 
\define\pr{\text{{\rm pr}}} 
\define\N{{\Bbb N}} 
\define\R{{\Bbb R}} 
\define\C{{\Bbb C}} 
\define\Z{{\Bbb Z}} 
\define\<{{\langle}} 
\define\>{{\rangle}} 
\define\minzero{\setminus\{0\}}

\define\gl{{g_\la}} 
\redefine\t{{\frak t}} 
\define\s{{\frak s}} 
\define\p{{\frak p}} 
\define\somr{{\SO(2m+2r)}} 
\define\sumr{{\SU(2m+r+1)}} 
\define\Inn{\text{\rm{Inn}}}

\define\sdp{{\medspace \times \kern -1.9pt 
         \vrule width0.4pt height 4.7pt 
          depth -0.3pt \medspace}} 
 
\define\restr#1{{\lower0.4ex\hbox{$\vert$}\lower0.7ex 
  \hbox{$\ssize{#1}$}}} 
\define\drestr#1{{\lower0.6ex\hbox{$\vert$}\lower1ex 
  \hbox{$\ssize{#1}$}}} 
\define\tdiffzero#1{{\frac d{d{#1}}\lower0.2ex\hbox{\restr{{#1}=0}}\,}} 
\define\diffzerosec#1{{\frac{d^2}{d{#1}^2}\lower1ex\hbox{$\big 
  \vert_{\raise1\jot\hbox{${\tsize {#1}=0}$}}$}\,}} 

\topmatter

\title Isospectral potentials and conformally equivalent isospectral
metrics on spheres, balls and Lie groups
\endtitle

\rightheadtext{Isospectral potentials and conformally equivalent
metrics}

\author Carolyn S. Gordon and Dorothee Schueth\endauthor

\address Dartmouth College, Hanover, New Hampshire 03755, U.S.A.
\endaddress

\email csgordon\@dartmouth.edu\endemail

\address Mathematisches Institut, Universit\"at Bonn, Beringstr\. 1,
         D@-53115 Bonn, Germany\endaddress

\email schueth\@math.uni-bonn.de\endemail

\thanks The first author is partially supported by National
  Science Foundation grant DMS-0072534.
  The second author is partially supported by SFB~611, Bonn.
\endthanks


\dedicatory
To the memory of our colleague and friend Robert Brooks (1952--2002),
whose influence is felt throughout our work.
\enddedicatory

\abstract
We construct pairs of conformally equivalent isospectral Riemannian
metrics $\varphi_1g$ and $\varphi_2g$ on spheres $S^n$ and balls $B^{n+1}$
for certain dimensions $n$, the smallest of which is $n=7$, and on
certain compact simple Lie
groups.  In the case of Lie groups, the metric~$g$ is left-invariant.
In the case of spheres and balls, the metric~$g$ is not the standard
metric but may be chosen arbitrarily close to the standard one.
For the same
manifolds $(M,g)$ we also show that the functions $\varphi_1$ and $\varphi_2$
are isospectral potentials for the Schr\"odinger operator
$\hbar^2\Delta +\varphi$.
To our knowledge, these are the first examples of isospectral potentials
and of isospectral conformally equivalent metrics on simply connected
closed manifolds.
\endabstract

\endtopmatter

\document

\heading Introduction\endheading

\noindent
Let $M$ be a compact Riemannian manifold.  Two fundamental questions in
spectral geometry are:

\item{$\bullet$}  To what extent can  one recover the geometry of $M$
from the eigenvalue spectrum of the Laplace-Beltrami operator
$\Delta$ acting on smooth functions on
$M$?  (If the manifold has boundary, one may impose either Dirichlet or
Neumann boundary conditions on the functions.)

\item{$\bullet$}  To what extent does the eigenvalue spectrum of the
Schr\"odinger operator $\hbar^2 \Delta + \phi$ determine the potential
$\phi\in C^\infty(M)$?

\noindent
Various geometric invariants of Riemannian manifolds such as dimension,
volume and total scalar curvature are known to be spectrally determined. 
The primary methods for identifying spectral invariants is through the
study of the heat and wave equations. If two Schr\"odinger operators
$\hbar^2\Delta+\phi$ and $\hbar^2\Delta+\psi$ on a fixed Riemannian
manifold~$M$ are isospectral, then the heat invariants show that $\phi$
and~$\psi$ have the same average and the same $L^2$@-norm (see,
e.g.,~\cite{Gi}).  

We will say two compact Riemannian manifolds are {\it
isospectral\/} if the associated Laplacians have the same eigenvalue
spectrum.  We will say that two potentials
$\phi_1$ and $\phi_2$ are {\it isospectral} if the Schr\"odinger
operators $\hbar^2 \Delta + \phi_1$ and $\hbar^2\Delta + \phi_2$ are
isospectral for every choice of the Planck's constant $\hbar$.  The only
way of identifying specific geometric invariants, respectively properties
of potentials, which are not spectrally determined is through construction
of isospectral metrics or potentials.  The primary results of this article
are the following constructions:
\roster
\item"(A)"
Isospectral conformally equivalent Riemannian metrics on
simply-connected manifolds, including spheres of dimension
seven and higher, balls of dimension at least eight, and on the Lie
groups $\Spin(2m+4)$, $m\geq 5$ and $\SU(2m+3)$, $m\geq 3$.  The metrics on
the spheres, respectively balls, may be chosen arbitrarily close to the
round, respectively flat metrics.  The metrics on the Lie groups are
conformally equivalent to left-invariant metrics.
\item"(B)"
Isospectral potentials for the Schr\"odinger operators
associated with (non-canonical) Riemannian metrics on the same underlying
simply-connected manifolds as in (A).  Again the Riemannian metrics on
the balls and spheres may be chosen arbitrarily close to the standard
ones, and the Riemannian metrics on the Lie groups are left-invariant.
\endroster
To our knowledge, these constructions give the first examples of
isospectral conformally equivalent metrics and of isospectral potentials
on simply-connected closed manifolds.  

We briefly review the literature on isospectral metrics and potentials. 
(We do not give a complete list of all known examples.)

The known methods for constructing isospectral manifolds roughly fall
into three classes:  ``Explicit'' computations, representation
theoretic or trace formula methods, and the use of torus actions. Among
the examples constructed by more or less explicit computations are 
$16$@-dimensional tori constructed by J. Milnor
\cite{Mi} in 1964 (these are the first known examples of isospectral
manifolds), spherical space forms \cite{Ik1,2}, domains in $\R^4$
\cite{Ur}, flat manifolds exhibiting various interesting properties
\cite{MR1-3}, and the first examples of isospectral manifolds with
different
local geometry
\cite{Sz1} and of pairs of isospectral metrics on spheres \cite{Sz2}. 
Representation methods were first used to construct isospectral Riemann
surfaces \cite{Vi} and continuous isospectral deformations \cite{GW1}.   
In 1985, T. Sunada \cite{Su} gave an elegant trace formula method
(which may also be
viewed as a representation theoretic method) for constructing isospectral
manifolds with a common finite cover.  This method and various
generalizations
\cite{Be1,2}, \cite{DG}, \cite{Pe1,2} led to a vast literature of
isospectral manifolds, including isospectral Riemann surfaces of every
genus greater than three
\cite{BT}, \cite{Bu}, \cite{BGG}, isospectral plane domains \cite{GWW},
and many other examples.
The third method, that of torus actions, was developed in the series of
papers \cite{Go1,3}, \cite{GW2}, \cite{GGSWW}, \cite{GSz},
\cite{Sch1-3}; see also \cite{Ba}.
This method generally produces manifolds with
different local geometry, including among others
the first examples of isospectral
simply-connected manifolds (in fact, continuous deformations) \cite{Sch1},
isospectral deformations of left-invariant
metrics on simply-connected compact Lie groups
\cite{Sch2},  and isospectral deformations of metrics on spheres
\cite{Go2}, \cite{Sch3}.
Recently, C.~Sutton \cite{St} modified Sunada's method (the second method
above) to construct new and very different examples of isospectral
simply connected manifolds.

R. Brooks, P. Perry and P. Yang \cite{BPY} gave a modification of Sunada's
method which yields isospectral conformally equivalent Riemannian metrics
with the same local but different global geometry. (This method is
closely related to an earlier modification of Sunada's method by R.
Brooks \cite{Br} which produces isospectral potentials for the
Schr\"odinger operator.)  The second author \cite{Sch2} analogously
extended the method of torus actions to construct the first examples of
isospectral conformally equivalent metrics that do not lift to isometric
metrics on the universal covering.  These were conformally equivalent to
left-invariant metrics on products
$G\times T$, where $G$ was a simple compact Lie group ($\Spin(n)$,
$\SO(n)$
or $\SU(n)$ for various $n$) and $T$ was a torus. These examples,
together with the
constructions of isospectral metrics on Lie groups and spheres
referred to above,  motivated our construction (A).

For the Schr\"odinger operator, there is a vast  literature concerning
isospectral potentials on the circle and on the interval. See, for example,
\cite{PT}.  An especially interesting connection
exists with the Korteveg de Vries equation.  Letting an arbitrary
function evolve according to the Korteveg de Vries equation, one obtains
a continuous  family of potentials $\phi_t$ on the circle such that the
operators
$\Delta+\phi_t$ are mutually isospectral
\cite{GGKM}. In higher dimensions, Brooks' modification of Sunada's
method, cited above, has yielded isospectral potentials on negatively
curved manifolds \cite{Br} and isospectral deformations of potentials
\cite{BG}. Although it was not mentioned explicitly in those two papers,
these
potentials are isospectral in the strong sense defined above; i.e., they
are isospectral for the Schr\"odinger operators with arbitrary Planck
constant $\hbar$.  In contrast to the potentials that we will construct,
Brooks' method yields potentials that lift to congruent  potentials on the
universal covering manifold.

The paper is organized as follows:
In Section~1 we present the general principles which will be
used for constructing
conformally equivalent isospectral metrics and isospectral potentials.
In Section~2 we outline a strategy for proving nonisometry.
Sections~3 and~4 contain the actual examples: Lie groups in Section~3,
balls and spheres in Section~4. In Section~5 we prove nonisometry
results for the isospectral metrics and potentials constructed in
the previous two sections. The stategy outlined in Section~2 is
applied in the case of the ball and
sphere examples; a different approach is used in the Lie groups case.
Finally, in Section~6 we give a short
survey of how our technique may be applied in other settings.

It is a pleasure to thank Steven Zelditch for helpful comments.

As we were completing this paper, we learned with shock and sorrow of the
untimely passing of our close friend and colleague, Robert Brooks.
Brooks brought to inverse spectral geometry ideas from such diverse areas
as complex analysis, probability and combinatorics.  His generosity and
enthusiasm in sharing his insights has profoundly affected our research.
He is deeply missed.

\bigskip

\heading\S1 Principles for constructing isospectral potentials and
conformally equivalent isospectral metrics\endheading

\definition{1.1 Definition}

(i) The {\it spectrum\/} of a closed Riemannian manifold $(M,g)$ is the
spectrum
of eigenvalues, counted with multiplicities, of the associated Laplace
operator~$\Delta_g$ acting on functions. The {\it Dirichlet spectrum\/}
of a compact
Riemannian manifold~$(M,g)$ with boundary is the spectrum of eigenvalues
corresponding to eigenfunctions which satisfy the Dirichlet boundary
condition $f\restr{\partial M}=0$.
The {\it Neumann spectrum} of such a manifold
is defined analogously with respect to the Neumann boundary condition
$Nf=0$, where~$N$ is the inward-pointing unit normal field on the boundary.
Given any potential function $\phi\in C^\infty(M)$ and a Planck constant
$\hbar\ne0$, one analogously
defines the (Dirichlet, Neumann) spectrum of the operator
$\hbar^2\Delta_g+\phi:f\mapsto
\hbar^2\Delta_gf+\phi f$.

(ii) Two closed Riemannian manifolds are called {\it isospectral\/}
if they have the same spectrum (including multiplicities). Two compact
Riemannian manifolds with boundary are called {\it Dirichlet isospectral},
resp\. {\it Neumann isospectral}, if they have the same Dirichlet spectrum,
resp\. the same Neumann spectrum. Isospectrality of two operators
$\hbar^2\Delta_g+\phi$ and $\hbar^2\Delta_{g'}+\psi$ is defined
analogously. If $(M,g)=(M',g')$ and if
$\hbar^2\Delta_g+\phi$ and $\hbar^2\Delta_g+\psi$ are (Dirichlet, Neumann)
isospectral for each $\hbar\ne0$, then
$\phi$ and~$\psi$ are called (Dirichlet, Neumann) isospectral potentials
on $(M,g)$.
\enddefinition

\subheading{1.2 Notation}
By a {\it torus}, we always mean a nontrivial, compact, connected
abelian Lie group. If a torus~$T$ acts smoothly and effectively by
isometries on a compact connected Riemannian manifold $(M,g)$, then
we denote by~$\hat M$ the union of those orbits on which~$T$ acts freely.
Note that~$\hat M$ is an open dense submanifold of~$M$.
The action of~$T$ gives~$\hat M$ the structure of a principal $T$@-bundle.
By~$g_0^T$ we denote the unique Riemannian metric on the quotient manifold
$\hat M/T$ such that the canonical projection $\pi:(\hat M,g_0)\to
(\hat M/T,g_0^T)$ is a Riemannian submersion.

\proclaim{1.3 Theorem}
Let $T$ be a torus which acts effectively by isometries on two compact
connected Riemannian manifolds $(M,g)$ and $(M',g')$. Let $\phi\in
C^\infty(M)$ and $\psi\in C^\infty(M')$ be two potential functions.
For each subtorus~$W\subset T$ of codimension one, suppose that there
exists a $T$@-equivariant diffeomorphism $F_W: M\to M'$ which
satisfies:
\roster
\item"(a)"
$F_W^*\dvol_{g'}=\dvol_{g}$\,,
\item"(b)"
$F_W$ induces an isometry $\bar F_W:(\hat M/W,g^W)\to
(\hat M'/W,g^{\prime\,W})$,
\item"(c)"
$F_W^*\psi=\phi$.
\endroster
Then for each $\hbar\ne0$, the operators $\hbar^2\Delta_g+\phi$ on~$M$ and
$\hbar^2\Delta_{g'}+\psi$ on~$M'$
are isospectral; if the manifolds have
boundary then these operators are Dirichlet and Neumann isospectral.

In particular, if moreover there exists an isometry $\tau:(M,g)\to(M',g')$
such that $\tau^*\psi\ne\phi$,
then $\phi$ and $\tau^*\psi$ are {\rm(}Dirichlet and Neumann{\rm)}
isospectral
potentials on~$(M,g)$.
\endproclaim

Note that for $\phi=\psi=0$, this is just Theorem~1.4 from~\cite{Sch3},
which in turn was  a slightly improved version of Theorem~1.2
from~\cite{Go3}.

\demo{Proof}
Let $\HH$ denote the Sobolev space $H^{1,2}(M,g)$ in the case of
closed~$M$ or in the case of Neumann boundary
conditions, and let ~$\HH$ denote the closure of
$C_0^\infty(M)$ in $H^{1,2}(M,g)$ in the case of Dirichlet boundary
conditions. For any subtorus~$W$ of~$T$ let $\HH_W\subset\HH$ denote the
subspace of $W$@-invariant functions. Analogously, we define the
function spaces~$\HH'$ and~$\HH'_W$ on~$M'$.

The key point in the proof of Theorem~1.4 from~\cite{Sch3}
(corresponding to
the case $\phi=\psi=0$) was deducing from~(a) and~(b) that for each~$W$
of codimension~$1$, the
pullback
$$F_W^*:\HH'_W\to\HH_W
$$
preserves both the $L^2$@- and the
$H^{1,2}$@-norm. Combining these maps according to the
Fourier decomposition
$$\HH=\HH_T\oplus\tsize\bigoplus_W(\HH_W\ominus\HH_T),
$$
one thus obtained an
$L^2$@-norm preserving Hilbert space isometry~$I:\HH'\to\HH$.
The statement (for $\phi=\psi=0$) then followed from the variational
characterization of eigenvalues via the Rayleigh quotient.

We need only slightly extend this proof in order to obtain Theorem~1.3:
By condition~(c), together with~(a), the pullbacks $F_W^*$ now also
satisfy $\int_M\phi|F_W^*f|^2\dvol_g$ $=$ $\int_{M'}\psi|f|^2\dvol_{g'}$
for all $f\in\HH'_W$\,. In particular, for the $L^2$@-norm preserving
isometry~$I:\HH'\to\HH$ defined above, we now have
$$\RR_\phi(If)=\RR_\psi(f)
$$
for all $f\in\HH'\minzero$, where
$$\align\RR_\psi(f):&=\tsize\int_{M'}(h^2\|df\|_{g'}^2+\psi|f|^2)
\dvol_{g'}\Big/
\tsize\int_{M'}|f|^2\dvol_{g'}\\
&=\bigl((h^2\|f\|^2_{H^{1,2}(M',g')}+\tsize\int_{M'}\psi|f|^2
\dvol_{g'})\bigm/\|f\|^2_{L^2(M',g')}\bigr)-h^2.\endalign
$$
is the modified Rayleigh quotient occurring in the variational
characterization of the eigenvalues of $\hbar^2\Delta_{g'}+\psi$, and
where
$\RR_\phi$ is defined analogously.

Thus the theorem follows again from the variational characterization
of eigenvalues.
\qed\enddemo

\subheading{1.4 Notation}
We fix a torus~$T$ with Lie algebra~$\z$.
Let~$\LL$ be the cocompact lattice in~$\z$ such that $\exp:\z\to T$
induces an isomorphism from $\z/\LL$ to~$T$, and denote by~$\LL^*
\subset\z^*$ the dual lattice. We also fix a compact connected Riemannian
manifold
$(M,g_0)$, with or without boundary, and a smooth effective action of~$T$
on $(M,g_0)$ by isometries.
\roster
\item"(i)"
For $Z\in\z$ we denote by $Z^*$ the vector field
$x\mapsto\tdiffzero t\,\exp(tZ)x$ on~$M$.
For each $x\in M$, we let
$\z_x:=\{Z^*_x\mid Z\in\z\}\subset T_xM$.
\item"(ii)"
We call a smooth  $1$@-form on~$M$ {\it admissible\/}
if it is $T$@-invariant and horizontal (i.e., vanishes
on the vertical spaces~$\z_x$).
\item"(iii)"
For any admissible $\z$@-valued $1$@-form~$\la$ on~$M$ we denote
by~$g\subla$ the Riemannian metric on~$M$ given by
$$g\subla(X,Y)=g_0\bigl(X+\la(X)^*,Y+\la(Y)^*\bigr).
$$
Note that $g\subla$ is the pullback of~$g_0$ by a field of unipotent
endomorphisms; in particular, $\dvol_{g\subla}=\dvol_{g_0}$\,.
Moreover, $g\subla$ is again $T$@-invariant.
\item"(iv)"
We say that a diffeomorphism $F:M\to M$ is {\it $T$@-preserving\/}
if conjugation by~$F$ preserves $T\subset\text{Diffeo(M)}$.
In that case, we denote by~$\Psi_F$ the automorphism of~$\z$
induced by conjugation by~$F$. Note that if $F$ is $T$@-equivariant
then $\Psi_F=\Id$. Obviously, each $T$@-preserving
diffeomorphism~$F$ of~$M$ maps $T$@-orbits to $T$@-orbits and satisfies
$F_*(Z^*)=\Psi_F(Z)^*$ for all $Z\in\z$, where the vector fields~$Z^*$
on~$M$ are defined as in~(i).

\endroster

\smallskip

The following theorem gives a useful specialization of Theorem~1.3
and is the generalized analog of Theorem~1.6 from~\cite{Sch3},
which corresponds to statement~(i) below in the special case $\phi=0$.

\proclaim{1.5 Theorem}
In the context of Notation~{\rm 1.4},
let $\la$, $\lap$ be two admissible $\z$@-valued $1$@-forms on~$M$.
Assume:
\roster
\item"($*$)" For every $\mu\in\LL^*$ there exists a $T$@-equivariant
  $F_\mu\in\Isom(M,g_0)$ such that $\mu\circ\la=F_\mu^*(\mu\circ\lap)$.
\endroster
Let $\phi\in C^\infty(M)$ be a potential function with the
property that the isometries~$F_\mu$ occurring in~\thetag{$*$} can be
chosen such that~$\phi$ is invariant under each of them.
Then we have:
\roster
\item"(i)"
$\hbar^2\Delta_{g\subla}+\phi$ and $\hbar^2\Delta_{g\sublap}+\phi$ are
isospectral for each $\hbar\ne0$;
if~$M$ has boundary then these operators are Dirichlet and Neumann
isospectral.
\item"(ii)"
If coincidentally there exists a $T$@-preserving $\tau\in\Isom(M,g_0)$
such that, in the notation of {\rm1.4(iv)},
$\Psi_\tau\circ\la=\tau^*\lap$, but
$\phi\ne\tau^*\phi$, then~$\phi$ and $\tau^*\phi$ are isospectral
potentials on $(M,g\subla)$; if $M$ has boundary then the two potentials
are Dirichlet and Neumann isospectral.
\endroster
\endproclaim

\demo{Proof}
(i) We claim that $(M,g\subla)$, $(M,g\sublap)$ and $\phi$, $\psi:=\phi$
satisfy the hypotheses of Theorem~1.3. Given a subtorus~$W$ of
codimension~$1$ in~$T$ we choose $\mu\in\LL^*$ such that $\kernn\mu=T_eW$,
and we choose a corresponding~$F_\mu$ as in~\thetag{$*$}.
In the proof of Theorem~1.6 from~\cite{Sch3} (corresponding to
the case $\phi=0$ here) it was shown that $F_W:=F_\mu$ satisfies
conditions~(a) and~(b) of Theorem~1.3; more precisely, (a) is immediate
from $\dvol_{g\subla}=\dvol_{g_0}=\dvol_{g\sublap}$\,, and (b) was deduced
from condition~\thetag{$*$}.
It only remains to check that~$F_\mu$ satisfies condition~(c) of
Theorem~1.3
(with $\psi=\phi$); but this is already guaranteed by the assumptions
on~$\phi$.

(ii) This follows immediately from~(i) by noting that~$\tau$ is an
isometry from $(M,g\subla)$ to $(M,g\sublap)$.
\qed\enddemo

\proclaim{1.6 Corollary}
Let $\la,\lap$ be two admissible $\z$-valued $1$@-forms on~$M$ satisfying
condition~\thetag{$*$} of Theorem~{\rm1.5}. Moreover, let
$\phi\in C^\infty(M,\R_+)$ be invariant under
the action of~$T$, and suppose that the isometries~$F_\mu$ occurring
in~\thetag{$*$} can be chosen such that~$\phi$ is invariant under
each of them. Then we have:
\roster
\item"(i)"
$(M,\phi g\subla)$ and $(M,\phi g\sublap)$ are isospectral.
\item"(ii)"
If coincidentally there exists a $T$@-preserving $\tau\in\Isom(M,g_0)$
such that \hbox{$\Psi_\tau\circ\la=\tau^*\lap$,} but $\tau^*\phi\ne\phi$,
then the conformally equivalent metrics $\phi g\subla$ and
$(\tau^*\phi)g\subla$ are isospectral.
\endroster
\endproclaim

\demo{Proof}
(i) Note that by definition, $\phi g\subla = (\phi g)\subla$
(and similarly for~$\lap$). Thus the statement follows by applying 
Theorem~1.5(i) (with trivial potential)
to $\la,\lap$ and to $\phi g_0$ instead of~$g_0$\,.

(ii) This follows immediately from~(i) by noting that
$\tau$ is an isometry from $(M,(\tau^*\phi)g\subla)$ to $(M,\phi g\sublap)$.
\qed\enddemo

\bigskip
\heading\S2 A nonisometry principle\endheading

\noindent
In the following we let $(M,g_0)$, $T$, $\z$ be as in Notation~1.4,
and we consider two admissible $\z$@-valued $1$@-forms~$\la,\lap$ on~$M$
in the sense of Notation~1.4(ii). We define~$\hat M$ as in Notation~1.2.

\subheading{2.1 Notation and Remarks} (\cite{Sch3})
\roster
\item"(i)"
We denote by $\Autm$ the group of all $T$@-preserving diffeomorphisms~$F$
of~$M$ which, in addition, preserve the $g_0$@-norm of vectors tangent
to the $T$@-orbits and induce an isometry, denoted~$\bar F$, of
$(\hat M/T,g_0^T)$. We denote the corresponding group of induced
isometries by $\Autbm\subset\Isom(\hat M/T,g_0^T)$.
\item"(ii)"
We define $\Cal D:=\{\Psi_F\mid F\in\Autm\}
\subset\Aut(\z)$ (see Notation~1.4(iv)).
Note that~$\Cal D$ is discrete because it is a subgroup of the discrete
group $\{\Psi\in\Aut(\z)\mid\Psi(\LL)=\LL\}$, where~$\LL$ is the lattice
$\kernn(\exp:\z\to T)$.
\item"(iii)"
Let~$\om_0$ denote the connection form on~$\hat M$ associated
with~$g_0$\,; i.e., for each $x\in\hat M$ the horizontal space
$\kernn(\om_0\restr{T_x\hat M})$ is the $g_0$@-orthogonal
complement of~$\z_x$ in~$T_x\hat M$ and $\om_0(Z^*_x)=Z$ for $Z\in\z$.
Then the connection form on~$\hat M$ associated with~$g\subla$ 
is given by $\om\subla:=\om_0+\la$.
Let~$\Om\subla$ denote the curvature form on $\hat M/T$ associated
with the connection form~$\om\subla$ on~$\hat M$.
We have $\pi^*\Om\subla=d\om\subla$ because $T$~is abelian.
\item"(iv)"
Since~$\la$ is $T$@-invariant and horizontal it induces a~$\z$@-valued
$1$@-form~$\bar\la$ on~$\hat M/T$. We conclude from $\pi^*\Om\subla=
d\om\subla=d\om_0+d\la$ that $\Om\subla=\Om_0+d\bar\la$. In particular,
$\Om\subla$~and~$\Om_0$ differ by an exact $\z$@-valued $2$@-form.
\endroster

\proclaim{2.2 Proposition (\cite{Sch3})}
Let~$\la,\lap$ be admissible $\z$@-valued $1$@-forms on~$M$ satisfying
the following two conditions:
$$\gather
\Om\subla\notin\Cal D\circ\Autbm^*\Om\sublap\tag{N}\\
\text{No nontrivial ${}\,1$@-parameter group in ${}\,\Autbm$
preserves~$\Om_{\lap}$\,.}\tag{G}
\endgather
$$
Then $(M,g\subla)$ and $(M,g\sublap)$ are not isometric.
\endproclaim
\subheading{Remark 2.3}
(i) For our examples of conformally equivalent isospectral metrics on
balls and spheres that we will construct in Sections 4 by applying
Proposition~1.6, we will use the nonisometry result above with
$\phi g_0$ playing the role of~$g_0$\,.  We will not be able to apply
this result to the Lie group examples in Section 3, however, and thus will
prove the non-isometry  of those examples by a different technique.

(ii) To construct isospectral potentials via Theorem~1.5, we may choose to
use the same functions~$\phi$ as for the examples of conformally
equivalent isospectral metrics arising from Corollary~1.6 (although
Theorem~1.5 permits a wider choice of potentials). Once we show that
a pair of conformally equivalent isospectral metrics
$\phi g\subla$ and $(\tau^*\phi)g\subla$ constructed as in
Corollary~1.6(ii) is nonisometric, it will follow immediately that the
isospectral (by Theorem~1.5) potentials $\phi$ and $\tau^*\phi$ for the
Schr\"odinger operator on $(M,g\subla)$ are not congruent under any
isometry of $(M,g\subla)$.  Indeed, if the potentials were congruent,
then the metrics
$\phi g\subla$ and
$(\tau^*\phi) g\subla$ would be isometric.
\smallskip

\bigskip

\heading\S3 Examples I: Lie groups\endheading

\noindent
The examples in this section are modelled after the construction of
isospectral left invariant Riemannian metrics on compact Lie groups given
in \cite{Sch2}.  

Let $\g$ be a Lie algebra.  Recall that $\ad_{\g}$
is a subalgebra of $\text{Der}(\g)$, the algebra of derivations of~$\g$.
The connected subgroup of $\Aut(\g)$ with Lie algebra $\ad(\g)$ is
denoted $\Inn(\g)$ and is called the group of inner automorphisms of~$\g$.
If $G$ is any connected Lie group with Lie algebra $\g$, then
$\Inn(\g)=\Ad_G$\,.  In particular, if $G$ is a matrix group and thus
$\g$ is a matrix algebra, then  $\Inn(\g)$ consists of all maps of~$\g$
given by conjugation by matrices in~$G$.  

Recall also that a Lie algebra is said to be compact if it is the Lie
algebra of a compact Lie group.

\definition{3.1 Definition}
Let  $\g$ be a compact 
Lie algebra, let $\z$ be a real vector space with an inner product, and
let
$j,j':\z\to\g$ be linear maps.

(i) $j$ and $j'$ will be said to be {\it isospectral\/} if
for each $Z\in\z$, there exists
$\alpha_Z\in \Inn(\g)$ such that
$j'(Z)=\alpha_Z(j(Z))$.  (We will focus on matrix algebras, in which case
this condition may be rewritten as $j'(Z)=\Ad(A_Z)(j(Z))=A_Zj(Z)A_Z^{-1}$
for some
$A_Z$ in the associated connected matrix group.)

(ii) $j$ and $j'$ will be said to be
{\it equivalent} if there exists
an automorphism~$\Phi$ of~$\g$ (not necessarily inner)
and an orthogonal linear map $C$ of $\z$ such that
$j'(Z)=\Phi(j(C(Z)))$ for all $Z\in\z$.  In this case we will say that
the pair $(\Phi,C)$ is an equivalence between $j$ and $j'$.
\enddefinition

\remark{Remarks} The notion of isospectrality and of equivalence of maps 
$j,j':\z\to\g$ was first introduced in \cite{GW2} in case $\g=\so(m)$
and in \cite{Sch2,3} in case $\g=\su(m)$.
Our notions of isospectrality and equivalence
in case $\g=\so(m)$ differ slightly from those in \cite{GW2}:
the conjugating elements $A_Z\in\SO(m)$ in
Definition~3.1(i) were allowed to lie instead in $\O(m)$ in \cite{GW2}.
In the notion of equivalence used in \cite{GW2}, the automorphism~$\Phi$
was required to lie in $\Ad_{\O(m)}$\,.  Since $\Aut(\so(m))=\Ad_{\O(m)}$
except when $m=4$ or $m=8$, the two notions of equivalence 
actually agree except in these special dimensions.
\endremark

\medskip

Isospectral maps $j$ and $j'$ will be used in various constructions of
isospectral Riemannian manifolds.  In general, equivalence  of the maps
$j$ and $j'$ will imply isometry of the resulting  manifolds; the
converse will hold under generic conditions.  

\proclaim{3.2 Proposition}
{\rm(i)} \cite{GW2} Let $m$ be any positive integer
other than $1,2,3,4$, or~$6$. Let $\Cal J$
be the real vector space consisting
of all linear maps $j:\R^2\to \so(m)$. Then there is a Zariski open
subset ${\Cal O}$ of~$\Cal J$ {\rm(}i.e., ${\Cal O}$ is
the complement of the zero
locus of some nonzero polynomial function on~$\Cal J${\rm)} such that
each $j \in {\Cal O}$ belongs to a $d$@-parameter family of isospectral,
inequivalent elements of~$\Cal J$.
Here $d\geq\hbox{$m(m-1)/2 - [m/2]([m/2]+2)$}>1$.
In particular, $d$~is of order at least $O(m^2)$.

{\rm(ii)} \cite{Sch2} Let $m\geq 3$. Let $\Cal J$ be the real vector
space consisting of all linear maps $j:\R^2\to \su(m)$. Then there is a
Zariski open subset ${\Cal O}$ of~$\Cal J$ such that each
$j \in {\Cal O}$ belongs to a continuous family of isospectral,
inequivalent elements of~$\Cal J$. 
\endproclaim

\remark{Remarks}
(i) The article \cite{GW2} also gives an example of a
continuous family of  isospectral, inequivalent maps from $\R^2$ to
$\so(6)$.

(ii) The two parts of the proposition were proven in \cite{GW2} and
\cite{Sch2,3}. In the case of $\so(m)$ in \cite{GW2}, the notions of
isospectrality and equivalence 
defined in the remark following Definition~3.1 were used.
However, the proposition remains
true if we instead define isospectrality as in
Definition~3.1. Indeed, as pointed out in
\cite{Sch2}, isospectrality (in the sense of
\cite{GW2}) of a continuous family $j_t$ as in Proposition~3.2(i)
means that on each
$Z\in\R^2$, the path $t\to j_t(Z)$ lies in a single $\Ad_{\O(m)}$-orbit
of~$j_0(Z)$ and thus, by continuity, the path in fact lives in a single
$\Ad_{\SO(m)}$-orbit.  Thus isospectrality also holds in the sense
of Definition~3.1.
A~similar argument may be used to show that inequivalence in the sense of 
\cite{GW2} of the continuous family $\{j_t\}$ implies inequivalence in the 
sense of Definition~3.1.

(iii) In Section 5, we will place various genericity conditions on the
maps $j$. The assertions of Proposition 3.2 remain true with these
genericity conditions, as one can see from the original proofs.
\endremark

\medskip

The following lemma is immediate.

\proclaim{3.3 Lemma}
Let $\k$ be a compact Lie algebra, and let $j_1,j_2:\z\to\k$ be
isospectral linear maps as in Definition~{\rm 3.1}.
Let $\g$ be a compact Lie algebra containing $\k\oplus\k$ as a
subalgebra.  Define $j, j':\z\to
\k\oplus\k\subset\g$ by
$j(Z)=(j_1(Z),j_2(Z))$ and $j'(Z)=(j_2(Z),j_1(Z))$. Then
$$j'(Z)=(\alpha_Z\,,\alpha_Z\inv)(j(Z))
$$
for each $Z\in\z$, where $\alpha_Z$ is given as in Definition~{\rm3.1} with
respect to the pair of isospectral maps $j_1\,,j_2$\,.  Thus $j$ and
$j'$ are  isospectral.
\endproclaim

\definition{3.4 Notation and Remarks}
Let $G$ be a compact Lie group with Lie algebra~$\g$,
and let $g_0$ be a bi-invariant Riemannian metric on~$G$.
Let $H<G$ be a closed connected subgroup of~$G$ with Lie algebra~$\h$,
and let $T<G$ be a torus with Lie algebra~$\z$.
Suppose that $[\z,\h]=0$, and that $\z$ is $g_0$@-orthogonal to~$\h$
in~$\z$.
\roster
\item"(i)" For any linear map $j:\z\to\h$ we define an associated
left invariant $\z$@-valued $1$@-form~$\la^j$ on~$G$ by letting
$g_0(\la^j(X),Z)=g_0(X,j(Z))$ for all $Z\in\z$ and $X\in\g$
(in other words, $\la^j:\g\to\z$ is the $g_0$@-transpose of
$j:\z\to\h\subset\g$).
\item"(ii)" We note that $\la^j$ is admissible in the sense of~1.4
with respect to the action of~$T$ on~$G$ from the right.
In fact, $\la^j$ vanishes on the $g_0$@-orthogonal complement
of~$\h$ and thus, in particular, on~$\z$.
The left invariant vector fields~$Z$
with $Z\in\z$ are just the vector fields~$Z^*$ (as in~1.4) induced
by the right action of~$T$. Moreover, $\la^j$ is invariant under this
action: Since $[\z,\h]=0$, the Lie algebra~$\h$ is $g_0$@-orthogonal
to $[\z,\g]$; therefore, $\la^j$ vanishes on $[\z,\g]$ and is hence
invariant under $\Ad_T$\,. Thus $\la^j$ is not only left invariant but
also right invariant under~$T$.
\item"(iii)" Since $\la^j$ is admissible with respect to the right
action of~$T$, we can define the associated metric~$g_{\la^j}$
on~$G$ as in Notation~1.4. Note that $g_{\la^j}$ is left invariant
(and right invariant under~$T$).
\endroster
\enddefinition

\proclaim{3.5 Lemma \cite{Sch2}}
In the notation of {\rm 3.1} and {\rm 3.4}, suppose that
$j,j':\z\to\h$ are isospectral. Let $\la,\lap:\g\to\z$ be the
associated left-invariant $1$@-forms.  Given a linear functional~$\mu$
on~$\z$, let $Z\in\z$ satisfy $\mu(W)=g_0(W,Z)$ for all $W\in \z$. 
Define
$F_\mu:=L_{A_Z}\circ R_{A_Z}\inv\in\Aut(G)$, where $A_Z\in H$ is such that
$j'(Z)=\Ad(A_Z)(j(Z))$ {\rm(}the existence of $A_Z$ follows from
Definition~{\rm 3.1)}. Then
$$\mu\circ\la=F_\mu^*(\mu\circ\lap);
$$
moreover, $F_\mu$ is $T$@-equivariant because $H$ and~$T$ commute,
and is a $g_0$@-isometry because $g_0$ is bi-invariant.
In particular, the associated Riemannian metrics $g\subla$ and~$g\sublap$
on~$G$ are isospectral by Theorem~{\rm 1.5(i)} \rm{(}with
$\phi=0${\rm)}.
\endproclaim

\proclaim{3.6 Theorem}
Let $G$ be a compact Lie group with Lie algebra~$\g$, let $H$ be
a compact Lie subgroup with Lie algebra of the form $\h=\k\oplus\k$ for
some Lie algebra $\k$, and let $T < G$ be a torus 
with Lie algebra~$\z$.  Suppose that $[\z,\h]=0$, and that $\z$ is
orthogonal to~$\h$ with respect to a bi-invariant Riemannian metric
$g_0$ on $G$.  Let
$j_1\,,j_2:\z\to\k$ be isospectral linear maps as in~{\rm3.1}. Define
$j:\z\to\k\oplus\k=\h$ by
$j(Z)=(j_1(Z),j_2(Z))$. Denote by~$g\subla$ the associated left invariant
metric on~$G$ defined as in~{\rm3.4}, where
$\la:=\la^j$\,. Let $\phi$ be a smooth function on~$G$ which is invariant
under conjugation by elements of~$H$.
Suppose that there exists an isometric automorphism~$\tau$ of
$(G,g_0)$ such that $\tau\restr T=\Id$ and such that $\tau_*$ restricts
to the map $(X,Y)\mapsto(Y,X)$ on $\k\oplus\k=\h\subset\g$.
Then:
\roster
\item"(i)"
$\phi$ and $\tau^*\phi$ are isospectral potentials on $(G,g\subla)$.
\item"(ii)"
If, in addition, $\phi$ is positive and right invariant under~$T$,
then $\phi g\subla$ and $(\tau^*\phi)g\subla$ are conformally equivalent
isospectral metrics on~$G$.
\endroster
\endproclaim

\demo{Proof}
Define $j':\z\to\k\oplus\k$ by $j'(Z):=(j_2(Z),j_1(Z))$. By Lemma~3.3,
$j$ and~$j'$ are isospectral. Let $\lap:=\la^{j'}$\,.
By Lemma~3.5, $\la$ and~$\lap$ satisfy condition~\thetag{$*$}
of Theorem~1.5 with respect to the right action of~$T$ on~$G$,
with each $F_\mu$ given as conjugation by an element of~$H$.
Since $\tau$ is a $T$@-equivariant $g_0$@-isometry and obviously
satisfies $\tau_*j(Z)=j'(Z)$ for all $Z\in\z$, if follows immediately
from the definition of $\la$ and~$\lap$ that $\tau^*\lap=\la$.
Thus statements (i) and~(ii) follow from Theorem~1.5(ii) and Corollary~1.6,
respectively.
\qed\enddemo

\subheading{3.7 Examples}
(i) Given a vector space $\z$, say of
dimension $r$, we may realize $\z$ as the Lie algebra of a maximal torus~$T$
in $\SO(2r)$.  Given isospectral maps $j_1\,,j_2:\z\to\so(m)=:\k$, let
$G=\SO(2m+2r)$ with a bi-invariant Riemannian metric~$g_0$\,.
Then $G$ contains $\SO(m)\times \SO(m)\times\SO(2r)$ as the subgroup
formed by three diagonal blocks, which in turn contains $\SO(m)\times
\SO(m)\times T$. Set $H=\SO(m)\times \SO(m)$.  The map $(p,q)\mapsto(q,p)$
of
$H$ extends to an inner automorphism --- hence an
isometry --- $\tau$ of~$G$ (namely, conjugation by the orthogonal
map of $\R^{2m+2r}$ interchanging the first $m$ coordinates with the
next $m$ coordinates).
Thus Theorem~3.6 yields an associated
left invariant metric~$g\subla$ on~$G$ and,
for any smooth function~$\phi$ on~$G$ which is invariant under
conjugation by elements
of $H$, a pair of isospectral potentials for the
Schr\"odinger operator on
$(G,g\subla)$.  Choosing $\phi$ also to be positive and
right invariant under~$T$, we further obtain a pair of isospectral metrics
conformally equivalent to~$g\subla$\,. In particular, Proposition~3.2
and the subsequent remark allow us to construct isospectral potentials
and conformally equivalent isospectral
metrics on $\SO(2m+4)$ for all $m\geq 5$.
(Although Proposition~3.2 yields continuous families of isospectral
$j$@-maps, the construction in Theorem~3.6 requires us to work with pairs
of isospectral maps. Thus we do not obtain continuous isospectral
deformations of potentials or of conformally equivalent metrics by this
construction.)

(ii)
Consider the two-fold covering $\tilde G=\Spin(2m+2r)$ of $\SO(2m+2r)$. 
Letting $\tilde H$ be the connected subgroup of $\tilde G$ with Lie
algebra $\h=\so(m)\oplus\so(m)$ and $\tilde T$ the toral subgroup with
Lie algebra $\z\subset\so(2r)$, and letting 
$j_1\,,j_2:\z\to\k=\so(m)$ be isospectral linear maps, we then obtain,
exactly as in Example~(i), an associated metric~$g\subla$ on~$\tilde G$
and, for any function~$\phi$ on~$\tilde G$ which is invariant under
conjugation by elements of $\tilde H$\,, a
pair of isospectral potentials for the Schr\"odinger operator on
$(\tilde G,g\subla)$. Note for that purpose that the automorphism~$\tau$ 
of~$SO(2m+2r)$
lifts to an automorphism~$\tilde\tau$ of~$\tilde G$. If $\phi$ is also
positive and right invariant under~$\tilde T$ then we further obtain
a pair of isospectral metrics conformally equivalent to~$g\subla$\,.
Thus we get pairs of isospectral potentials and conformally equivalent
isospectral metrics on $\Spin(2m+4)$ for all $m\ge5$.

(iii) In the notation of Example (i), we may also realize~$\z$ as a
maximal torus in $\SU(r+1)$.  The analogous construction with $\k=\su(m)$
then yields, for
each pair of isospectral maps $j_1\,,j_2:\z\to\su(m)=:\k$, an associated
left invariant Riemannian metric $g_\la$ on $G:=\SU(2m+r+1)$ and,
for appropriately chosen~$\phi$, isospectral potentials on $(G,g\subla)$
and isospectral metrics on~$G$ conformally equivalent to $g\subla$\,.
In particular, applying Proposition~3.2(ii), we obtain such isospectral
potentials and conformally equivalent metrics on $\SU(2m+3)$ for all
$m\geq 3$.

\bigskip

\subheading{3.8 Remark}
In 5.1 we will prove nontriviality of the above examples from~3.7 for a
particular choice of $\phi$ under the hypothesis that $j_1$ is not
equivalent to $j_2$ and under certain genericity conditions on $j$.  

\bigskip

\heading\S4 Examples II: Spheres and Balls\endheading

\noindent
Let $r\ge2$, and let $T$ be the $r$@-dimensional torus $\R^r/\LL$
with $\LL=(2\pi\Z)^r$.
Let $\{Z_1\,,\ldots\,,Z_r\}$ be the standard basis of its Lie algebra
$\z=\R^r$.
\subheading{4.1 Notation}

(i) For any given $n$, we consider the action of~$T=(S^1)^r$ on
$\R^n\oplus\C^r$ which is induced by the canonical action of
$S^1\times\ldots
\times S^1$ on $\C\times\ldots\times\C$ given by
multiplication in each
factor.

(ii) Given any $\z$@-valued $1$@-form~$\la$ on~$\R^n$, its trivial
extension
to $\R^n\oplus\C^r$ (again denoted~$\la$) is admissible in the sense
of Notation~1.4 with respect
to the action of~$T$ defined above. Let $(B,g_0)$, respectively $(S,g_0)$,
be the standard unit ball, respectively sphere, in $\R^n\oplus\C^r\cong
\R^{n+2r}$. We continue to denote by~$\la$ the restriction, respectively
pullback, of~$\la$ to~$B$, respectively~$S$. Then~$\la$ is admissible with
respect to the induced action of~$T$ and thus defines a metric~$g\subla$
on~$B$ and on~$S$.

\medskip

\definition{4.2 Definition}
Two $\z$@-valued $1$@-forms~$\nu,\nup$ on~$\R^m$ are called
{\it isospectral\/} if for each $\mu\in\LL^*$ there exists $A_\mu\in\O(m)$
such that $\mu\circ\nu=A_\mu^*(\mu\circ\nu)$.
\enddefinition

\proclaim{4.3 Theorem}
Let $n=2m$, and
let $M$ denote either the ball~$B$ or the sphere~$S$ in
$\R^{n+2r}\cong\R^n\oplus\C^r$.
Let $\nu,\nu'$ be isospectral $\z$@-valued $1$@-forms on~$\R^m$.
Define a $\z$@-valued $1$@-form~$\la$ on $\R^n=\R^m\oplus\R^m$ by
$$\la_{(p,q)}(X,Y):=\nu_p(X)+\nu'_q(Y).
$$
Denote by~$g\subla$ the associated metric on~$M$ as in
Notation~{\rm4.1(ii)}.
Let $\psi:\R^m\to\R$ be a smooth radial function and define
$\phi_1\,,\phi_2:\R^m\oplus\R^m\oplus\C^r\to\R$ by
$\phi_i(p_1\,,p_2\,,u):=\psi(p_i)$ for $i=1,2$.
Then
\roster
\item"(i)"  $\phi_1$ and $\phi_2$ are  isospectral  potentials on 
$(M,g_\la)$.  
\item"(ii)" If $\psi$ is strictly positive, then $\phi_1g_\la$ and 
$\phi_2g_\la$ are isospectral conformally equivalent Riemannian metrics 
on~$M$.
\endroster
\endproclaim

\demo{Proof}
We apply the criteria of Theorem~1.5 and Corollary~1.6.
Define a second $\z$@-valued $1$@-form~$\lap$ on $\R^n=\R^m\oplus\R^m$
by
$$\lap_{(p,q)}(X,Y):=\nu'_p(X)+\nu_q(Y).
$$
Let $g\sublap$ be the associated Riemannian metric on~$M$.
Given any $\mu\in\LL^*$, choose $A_\mu\in\O(m)$ as in Definition~4.2.
Then the map
$$F_\mu:\R^m\oplus\R^m\oplus\C^r\ni (p,q,u)\mapsto(A_\mu p,A_\mu\inv q,u)
\in\R^m\oplus\R^m\oplus\C^r
$$
restricts to a $T$@-equivariant isometry of $(M,g_0)$ satisfying
$\mu\circ\la=F_\mu^*(\mu\circ\lap)$.
Moreover, the map
$$\tau:\R^m\oplus\R^m\oplus\C^r\ni(p,q,u)\mapsto(q,p,u)
\in\R^m\oplus\R^m\oplus\C^r
$$
also restricts to a $T$@-equivariant isometry of $(M,g_0)$ and satisfies
$\la=\tau^*\lap$.

We let $\phi_1$ play the role of~$\phi$ in the notation of~1.5 and~1.6,
and observe that $\phi_2=\tau^*\phi_1$\,.
Since $\psi$ is radial, $\phi_1$
is invariant under each~$F_\mu$\,. Thus Theorem~3.8 follows from
Theorem~1.5 and Corollary~1.6.
\qed\enddemo

\subheading{4.4 Example}
\noindent
Our first family of examples is modelled after the construction
of isospectral balls and spheres in~\cite{Go3} (see Remark~4.5 below).

Let $j,j':\z\to\so(m)$ be isospectral linear maps, where we allow the
weaker notion of isospectrality given in the remark following
Definition~3.1. Define two
$\z$@-valued
$1$@-forms
$\nu,\nup$ on~$\R^m$ by
$$\<\nu^{(\prime)}(X_p),Z\>=\<j^{(\prime)}(Z)p,X\>
$$
for  $X_p\in T_p\R^m$ and $Z\in\z$, where $\scp$ denotes the standard
inner products both on~$\R^m$ and on $\z=\R^k$.
Then $\nu$ and~$\nup$ are isospectral in the sense of Definition~4.2.
In fact, let~$Z\in\z$ be the vector dual to~$\mu$ with  respect
to the inner product $\scp$ on $\z$. Choose $A_Z\in\O(m)$ as in the
remark following Definition~3.1; then $A_\mu:=A_Z$ will satisfy
$\mu\circ\nu=A_\mu^*(\mu\circ\nup)$.

Recall from Proposition~3.2 that for $r=2$, i.e. $\z=\R^2$, and each
$m\ge5$ there do exist pairs $j,j':\z\to\so(m)$ as above.
Using the corresponding forms $\nu,\nup$ for constructing~$\la$ in
Theorem~4.3, we obtain a metric~$g\subla$\,,
pairs of isospectral metrics conformally
equivalent to~$g\subla$\,,
and pairs of potentials which are isospectral with respect to~$g\subla$
on each of the spheres $S^{2m+3\ge13}$ and the balls $B^{2m+4\ge14}$.

\smallskip

\remark{4.5 Remark}
In 5.2 we will show that the examples from~4.4 with $\z=\R^2$ are
nontrivial provided that
\roster
\item"(i)" the maps $j,j'$ are not equivalent in the sense of
Definition~3.1,
\item"(ii)" $j$ and $j'$ are {\it generic\/} in the sense
that no nonzero element of~$\so(m)$ commutes with both $j(Z_1),j(Z_2)$,
or with both $j'(Z_1),j'(Z_2)$, and
\item"(iii)" the restriction to $\overline{B_1(0)}\subset\R^m$
of the smooth radial function~$\psi:\R^m\to\R_+$ attains
its maximum precisely in $p=0$.
\endroster
\endremark

\smallskip

\remark{Remark}
In the notation of Theorem~4.3 and Example 4.4 above,
the isospectral spheres and balls from \cite{Go3} were just the unit
spheres,
respectively balls, in $\R^m\oplus\C^r$, endowed with the associated
isospectral metrics $g_\nu$ and $g_{\nup}$\,.
\endremark

\medskip

Note that in Example~4.4, we had $r\ge2$ and $m\ge5$, where $r$ and $m$
are the dimensions used in Theorem~4.3. We now give an example with $r=2$
and $m=3$, yielding conformally equivalent isospectral metrics and
isospectral potentials on $S^9\subset\R^3\oplus\R^3\oplus\C^2$ and
on~$B^{10}$.
The approach is related to certain constructions from \cite{Sch2,3}.

\subheading{4.6 Example}

\noindent
Let $r=2$, thus $\z=\R^2$.  Let
$\Sym_0(\R^3)$ denote the space of symmetric traceless
real $3\times3\,$@-matrices.
For each linear map $c:\z\cong\R^2\to\Sym_0(\R^3)$ we define a
$\z$@-valued
$1$@-form $\nu$ on~$\R^3$ by letting
$$\<\nu(X_p),Z\>=\<c(Z)p\times p,X\>
$$
for  all $X\in T_p\R^3$ and $Z\in\z$,
where $\scp$ denotes the standard euclidean inner products on both~$\R^3$
and on~$\z$, and~$\times$ denotes the vector product in~$\R^3$.

The key observation is that there exist pairs of linear maps
$c,c':\z\to\Sym_0(\R^3)$ such that
\roster
\item"1.)"
$c$ and~$c'$ are {\it isospectral}, that is,
for each $Z\in\z$ there exists $E_Z\in\SO(3)$ such that
$c'(Z)=E_Zc(Z)E_Z\inv$\,;
\item"2.)"
$c$ and~$c'$ are {\it not equivalent}, where being equivalent means:
There exists $E\in\O(3)$
and $C\in\O(2)$ such that $c'(Z)=Ec(C(Z))E\inv$ for all $Z\in\z$;
\item"3.)"
$c$ and $c'$ are {\it generic\/}, that is, no nonzero element
of~$\so(3)$ commutes with both $c(Z_1)$ and $c(Z_2)$, and the analogous
statement holds for $c'$. 
\endroster

An explicit example of a pair $c,c'$ satisfying the conditions above
is given by 
$$
c(Z_1)=c'(Z_1)=\left(\smallmatrix -1&0&0\\0&0&0\\0&0&1
\endsmallmatrix\right),\ \
c(Z_2)=\left(\smallmatrix 0&1/\sqrt2&0\\1/\sqrt2&0&1/\sqrt2\\
0&1/\sqrt2&0\endsmallmatrix\right),\ \
c'(Z_2)=\left(\smallmatrix 0&0&1\\0&0&0\\1&0&0\endsmallmatrix
\right)
$$
(compare~\cite{Sch3}, Proposition~3.3.6; actually $c(Z_2)$ and~$c'(Z_2)$
were rescaled by~$\sqrt2$ there, but this does not affect any of the
three properties above).

It is easy to check that condition~1.)~implies isospectrality (in the
sense of~4.2) of the associated $\z$@-valued $1$@-forms
$\nu,\nu'$ on~$\R^3$.
Using these in Theorem~4.3, we obtain
pairs of isospectral potentials and pairs
of conformally equivalent isospectral
metrics on $S=S^9\subset\R^3\oplus\R^3\oplus\C^2$, resp\. on
$B=B^{10}$.

In 5.2 we will prove the corresponding nontriviality statements,
provided that the restriction to $\overline{B_1(0)}\subset\R^3$ of
the smooth radial function~$\psi:\R^3\to\R_+$ attains its maximum
precisely in $p=0$.

\medskip

Having thus obtained examples of conformally equivalent isospectral
metrics and isospectral potentials on $S^9$ and $B^{10}$,
we do a final twist in order to obtain such
examples even on $S^7$ and $B^8$.
The idea is to lift the forms $\nu,\nup$ used in Example~4.6 from~$\R^3$
to~$\C^2$ via the Hopf projection. The new $\la=\nu\oplus\nup$ will live
on $\C^2\oplus\C^2\cong\R^8$ instead of $\R^3\oplus\R^3$. In return,
the third component
of our previous ambient space $\R^3\oplus\R^3\oplus\C^2$ will become
unnecessary because we will
use a different torus action, living on $\C^2\oplus\C^2$ itself.
This construction is related to a construction from~\cite{Sch3}
of isospectral metrics on~$S^5\subset\C^2\oplus\C$ and on~$B^6$.

First we need some preparations.

\subheading{4.7 Notation}
(i) Let $r=2$, thus $\z=\R^2$ and $T=\R^2/(2\pi\Z)^2$.
For any $m\in\N$ we identify the real vector spaces $\C^m$
and~$\R^{2m}$ via the linear isomorphism which sends
$\{e_1\,,ie_1\,,\,\ldots,e_m\,,ie_m\}$ (in this order) to
the standard basis of~$\R^{2m}$, where $\{e_1\,,\,\ldots,e_m\}$
denotes the standard basis of~$\C^m$.
We let the torus~$T$ act on $\C^m\oplus\C^m$ by
$$\exp(aZ_1+bZ_2):(p,q)\mapsto(e^{ia}p,e^{ib}q)
$$
for all $a,b\in\R$, $p,q\in\C^m$.

(ii) Let $(B,g_0)$, resp.~$(S,g_0)$, be the standard unit ball
of dimension~$4m$, resp.~the standard unit sphere of dimension~$4m-1$,
in $\R^{4m}\cong\C^m\oplus\C^m$.
Given any $\z$-valued $1$@-form~$\la$ on $\C^m\oplus\C^m$ which is
admissible with respect to the action of~$T$ defined above, we also
denote by~$\la$ the restriction, resp.~the pullback, to~$B$, resp.~$S$.
Then $\la$ is admissible with respect to the induced action of~$T$ on~$B$
and~$S$, and thus defines a corresponding metric~$g_\la$ on~$B$,
resp.~$S$.

\smallskip

\definition{4.8 Definition}

(i)
A $\z$@-valued $1$@-form~$\nu$ on~$\C^m\cong\R^{2m}$ is called
{\it Hopf admissible\/} if it is admissible, in the sense of~1.4,
with respect to the Hopf action of the $1$@-dimensional Torus~$S^1$
on~$\C^m$.

(ii)
Two Hopf admissible $\z$@-valued $1$@-forms $\nu,\nu'$
on~$\C^m$ are called
{\it isospectral\/} if for each $\mu\in\LL^*$ there exists
$A_\mu\in\SU(m)$ such that $\mu\circ\nu=A_\mu^*(\mu\circ\nu')$.
\enddefinition

\proclaim{4.9 Theorem}
Let $M$ denote either the ball~$B$ or the sphere~$S$ in
$\C^m\oplus\C^m\cong\R^{4m}$.
Let $\nu$ and $\nu'$
be  Hopf admissible $\z$-valued $1$@-forms on~$\C^m$,
isospectral in the sense of Definition~{\rm4.8}.
Write $\nu=(\nu_1\,,\nu_2)$
and $\nu'=(\nu'_1\,,\nu'_2)$ with respect to the coordinates on $\z=\R^2$
defined by the basis $\{Z_1,Z_2\}$.  Assume that there exist
$A,A'\in\SU(m)$ such that $A^*\nu_1=\nu_2$\,, $A^*\nu_2=\nu_1$\,, and
similarly for~$A'$ with respect to~$\nu'$.
Define a $\z$@-valued $1$@-form~$\la$ on $\C^m\oplus\C^m$ by
$$\la_{(p,q)}(X,Y):=\nu_p(X)+\nu'_q(Y).
$$
for all $(p,q)\in\C^m\oplus\C^m$.
Let $\psi:\C^m\cong\R^{2m}\to\R$ be a smooth radial function and define
$\phi_1\,,\phi_2:\C^m\oplus\C^m\to\R$ by
$\phi_i(p_1\,,p_2):=\psi(p_i)$ for $i=1,2$.
Then
\roster
\item"(i)"  $\phi_1$ and $\phi_2$ are  isospectral  potentials on 
$(M,g_\la)$.  
\item"(ii)" If $\psi$ is strictly positive, then $\phi_1g_\la$ and 
$\phi_2g_\la$ are isospectral conformally equivalent Riemannian metrics 
on~$M$.
\endroster
\endproclaim

\demo{Proof}
The proof is very similar to that of Theorem~4.3, slightly complicated
by the fact that switching the components of $\C^m\oplus\C^m$ is
no longer a $T$@-equivariant map here.
Define a second $\z$@-valued $1$@-form~$\lap$ on $\C^m\oplus\C^m$ by
$$\lap_{(p,q)}(X,Y):=\nu'_p(X)+\nu_q(Y).
$$
Let $g\sublap$ be the associated Riemannian metric on~$M$.
Given any $\mu\in\LL^*$ and the associated $A_\mu$ as in
Definition~4.8(ii), the map
$$F_\mu:\C^m\oplus\C^m\ni (p,q)\mapsto(A_\mu p,A_\mu\inv q)
  \in\C^m\oplus\C^m
$$
restricts to an isometry of $(M,g_0)$ satisfying $\mu\circ\la=F_\mu^*
(\mu\circ\lap)$;
note that $F_\mu$ is $T$@-equivariant since each $A_\mu\in\SU(m)$
commutes with multiplication by complex scalars.
Moreover, the map
$$\tau:\C^m\oplus\C^m\ni(p,q)\mapsto(A'q,Ap)\in\C^m\oplus\C^m
$$
restricts to a $T$@-preserving isometry of $(M,g_0)$.
The corresponding
linear automorphism~$\Psi_\tau$ of~$\z$ interchanges the basis vectors
$Z_1$ and~$Z_2$\,. By the conditions on~$A$ and~$A'$, it is easy to
check that~$\tau$ satisfies $\Psi_\tau\circ\la=\tau^*\lap$.

We let $\phi_1$ play the role of~$\phi$ in the notation of~1.5 and~1.6,
and observe that $\phi_2=\tau^*\phi_1$ because~$\psi$ is radial.
By the same reason, $\phi_1$
is invariant under each~$F_\mu$\,. Thus Theorem~4.9 follows from
Theorem~1.5 and Corollary~1.6.
\qed\enddemo

\subheading{4.10 Example}

\noindent
We fix a realization of the Hopf projection $P:S^3\to S^2_{1/2}\subset
\R^3$ in coordinates, say
$$P:(\alpha,\beta,\gamma,\delta)\mapsto\bigl(\tfrac12(\alpha^2+\beta^2
  -\gamma^2-\delta^2),\,\alpha\gamma+\beta\delta,\,\alpha\delta-\beta
  \gamma\bigr).
$$
We extend~$P$ to a smooth map from~$\R^4\cong\C^2$ to~$\R^3$, defined
by the same formula (note that~$P$ will map~$S^3_a$ to~$S^2_{a^2/2}$
for each radius $a\ge0$).

With each linear map $c:\z\cong\R^2\to\Sym_0(\R^3)$ we associate a
$\z$@-valued $1$@-form $\nu=(\nu_1\,,\nu_2)$ on~$\C^2$ by letting
$$\nu_k(X_p)=\<c(Z_k)P(p)\times P(p),P_{*p}(X)\>
$$
for $k=1,2$ and all $X\in T_p\C^2$.

Consider the specific pair of linear maps $c,c':\R^2\to\Sym_0(\R^3)$
from Example~4.6. In addition to the three properties (isospectrality,
inequivalence, genericity) stated there, we observe that
\roster
\item"(I)"
There exist $E,E'\in\SO(3)$ satisfying $Ec(Z_1)E\inv=c(Z_2)$\,,
$Ec(Z_2)E\inv=c(Z_1)$\,, and similarly for~$E'$ with respect
to~$c'$.
\endroster
In fact, the elements $E,E'$ of $\SO(3)$
given by
$$
E=\left(\smallmatrix 1/2&-1/\sqrt2&1/2\\-1/\sqrt2&0&1/\sqrt2\\
1/2&1/\sqrt2&1/2\endsmallmatrix\right),\qquad
E'=\left(\smallmatrix -1/\sqrt2&0&1/\sqrt2\\0&-1&0\\1/\sqrt2&0&1/\sqrt2
\endsmallmatrix\right)
$$
satisfy $E^2=E^{\prime\,2}=\Id$ and conjugate $c(Z_1)$ to~$c(Z_2)$\,,
respectively $c'(Z_1)$ to~$c'(Z_2)$\,.

Again, condition~1.)~on $c,c'$ from Example~4.6 implies isospectrality
of the associated Hopf admissible $\R^2$@-valued $1$@-forms $\nu,\nup$
on~$\C^2$. In fact, for any given $\mu\in\LL^*\subset\z^*$
let $Z\in\z$ be the dual
vector (via the chosen basis). Choose $E_Z$ as in condition~1.),
and choose $A_Z\in\SU(2)\subset\SO(4)$ such that $P\circ A_Z=E_Z\circ P$.
Then it is straightforward to check that $A_\mu:=A_Z$ satisfies
$\mu\circ\nup=A_\mu^*(\mu\circ\nu)$.

Moreover, condition~(I) implies existence of elements $A,A'\in\SU(2)$
as required in the conditions of Theorem~4.9: Just choose $A,A'$ such that
$P\circ A=E\circ P$, $P\circ A'=E'\circ P$. (Note that $A,A'$ are both
unique up to sign and will satisfy $A^2=A^{\prime\,2}=-\Id=e^{i\pi}\Id$.)

Thus we are able to apply Theorem~4.9 and obtain
pairs of isospectral potentials and pairs
of conformally equivalent isospectral
metrics on $S=S^7\subset\C^2\oplus\C^2$ and on
$B=B^8$.
In~5.3 we will show nontriviality of these examples for every choice
of a positive smooth radial function $\psi:\R^4\cong\C^2\to\R_+$ which is
strictly monotonously increasing with the radius.

\subheading{4.11 Remark}
The choice of the metrics~$g\subla$
to which our isospectral metrics are conformal, resp\. which admit
the pairs of isospectral potentials constructed above, allows an
interesting
modification.
In each of our examples of this section (4.4, 4.6, 4.10)
one can actually {\it modify~$g\subla$ in such a way that it becomes
equal to the standard metric~$g_0$ outside certain subsets of
arbitrarily small volume.}
The idea, first used in~\cite{Sch3} for a similar aim, is to
exploit the following obvious fact:
\roster
\item"" If $\la,\lap$ are $\z$-valued admissible forms on~$M$ which
satisfy condition~\thetag{$*$} of Theorem~1.5, then so do $f\la,f\lap$,
where
$f\in C^\infty(M)$ is any function which is invariant under~$T$ and
under each~$F_\mu$ occurring in~\thetag{$*$}.
\endroster
In particular, the statements~1.5(i) and~1.6(i) will still hold with
$\la,\lap$ replaced by $f\la,f\lap$, and so will~1.5(ii) and~1.6(ii),
provided that~$f$ is invariant under~$\tau$, too.
In Theorems~4.3 and~4.9, which were specializations of~1.5/1.6,
we can therefore, without affecting any of the isospectrality statements,
replace~$\la$ by~$f\la$, where $f$ is the restriction to~$M$ of
\roster
\item"--" (for Theorem 4.3, hence the Examples~4.4 and~4.6:)
  any smooth function on $\R^m\oplus\R^m
  \oplus\C^r$ which is of the form $(p,q,u)\mapsto\alpha(|p|,|q|,|u|)$
  with $\alpha(r,s,t)=\alpha(s,r,t)$ for all $(r,s,t)$;
\item"--" (for Theorem 4.9, hence Example~4.10:)
  any smooth function on $\C^m\oplus\C^m$ which is of the form
  $(p,q)\mapsto\alpha(|p|,|q|)$ with $\alpha(a,b)=\alpha(b,a)$ for
  all~$(a,b)$.
\endroster

In particular, we may in each case choose~$\alpha$ to have as small
support (in the parameter set relevant for~$M$) as we like.
Then the support of~$f$ in~$M$ will have small volume with respect
to the standard metric~$g_0$ and hence also with respect to $g_{f\la}$
(recall from 1.4(iii) that $\dvol_{g_0}=\dvol_{g_{f\la}}$).

In 5.4, we will show that the corresponding nonisometry statements
will remain true, provided that $f$ does not vanish identically
on~$M$; only in the case
of Example 4.10 and $M=B^8$, we will need an additional technical
assumption: Namely,
that the support of~$\alpha$ nontrivially intersects a certain open
set $\Cal U\subset\{(a,b)\in B_1(0)\mid a>b>0\}$ (whose definition
will actually depend on the radial function~$\psi$ used for constructing
the potentials).

\bigskip

\heading\S5 Non-triviality of the examples\endheading
 
\subheading{5.1 Nonisometry proof for the examples from 3.7} 
 
\noindent 
We consider the construction of metrics and potentials on Lie groups 
given in Example~3.7(i--iii). 
We first prove that for a special choice of the function 
$\phi$ and with a very 
mild genericity condition on $j$ (condition 5.1.2(a) below), the 
isospectral  
potentials $\phi$ and $\psi:=\tau^*\phi$ on $(G,\gl)$ are 
not congruent if 
$j_1$ and $j_2$ are inequivalent (Theorem~5.1.8).  
We then show under an additional 
genericity hypothesis on $j=(j_1,j_2)$ (condition 5.1.2(b)) that the 
isospectral conformally equivalent metrics $\phi\gl$ and $\psi\gl$ are not 
isometric (Theorem~5.1.13). 
For notational simplicity, we will focus throughout on the two cases 
$G=\SO(2m+2r)$ and $G=\SU(2m+r+1)$.  With minor adjustments to the proofs 
in the case of $\SO(2m+2r)$, one obtains the same results for 
$G=\Spin(2m+2r)$. 
\definition{5.1.1 The function $\boldsymbol\varphi$} 
 
\noindent 
For any matrix $A\in M(m,m,\C)$ and any column vector $V\in\C^m$ we write 
$$d(A,V):=\det(V,AV,A^2 V,\ldots,A^{m-1} V). 
$$ 
Note that $d$ is invariant under the action of $\SL_\C(m)$ on 
$M(m,m,\C)\times\C^m$ given by $B(A,V)=(BAB\inv,BV)$ for $B\in\SL_\C(m)$. 
It is easy to see that if $A$ has pairwise different eigenvalues 
$\lambda_1\,,\ldots,\lambda_m$\,, then $d(A,V)$ is the Vandermonde 
determinant $\left|\Bigl(\;\lambda_k^\ell\;\Bigr)\right|_{ 
{1\le k\le m}\atop{0\le\ell\le m-1}} = \prod_{i<j}(\lambda_j-\lambda_i)$, 
multiplied by the product of the coefficients of~$V$ with respect to 
an eigenbasis of determinant one of~$A$. If $A$ is diagonizable and has 
at least one multiple eigenvalue, then $d(A,V)=0$ for all~$V$. 
 
For $A\in M(m,m,\C)$ and $C\in M(m,\ell,\C)$ we let 
$$d(A,C):=\prod_{k=1,\ldots,\ell} d(A,C_k), 
$$ 
where $C_1\,,\ldots,C_\ell$ denote the column vectors of~$C$. 
 
We write matrices in 
$G=\SO(2m+2r)$ or $G=\SU(2m+r+1)$ as 
$$X=\left(\matrix A&B&C\\D&E&F\\H&J&L\endmatrix\right), 
$$ 
where $A,B,D,E$ are $m\times m\,$@-matrices and $L$ is a 
$2r\times 2r\,$@-matrix or an $(r+1)\times(r+1)$@-matrix, respectively. 
Let $c_1>c_2>0$. In case $G=\somr$, we define~$\phi:G\to\R_+$ by 
$$\phi(X):=\exp(c_1\det A+c_2\det E). 
$$ 
In case $G=\sumr$, we let 
$$\phi(X):=\exp\bigl(c_1\,\text{\rm Re}\det A+c_2\,\text{\rm Re}\det E 
         -(\text{\rm Re}\,d(A,C))^2-(\text{\rm Re}\,d(E,F))^2\bigr), 
$$ 
where $\text{\rm Re}$ denotes the real part. 
\enddefinition 
 
\definition{5.1.2 Genericity conditions} 
\roster 
\item"(a)" 
The kernel of the  map $j=(j_1,j_2):\z\to\h$ is trivial, where 
$\h=\k\oplus\k$. 
\item"(b)" 
The  image of $j$ has trivial centralizer in $\h$. 
\endroster 
\enddefinition 
 
\definition{5.1.3 Notation}  
For any subgroup~$U$ of $G$, denote by $U_L$\,, respectively $U_R$\,, 
the group of all left, respectively right, translations of $G$ by elements 
of~$U$. Moreover, denote by $\Inn_G(U)$ the subgroup of the group of 
inner automorphisms of~$G$ consisting of conjugations by elements of~$U$. 
\enddefinition 
 
\proclaim{5.1.4 Lemma \cite{OT}} 
Let $G$ be a connected simple compact Lie group 
with a left invariant Riemannian metric~$g$.  Then: 
\roster 
\item"(i)" 
The full group of isometries  of $(G,g)$ is the semi-direct product of 
$G_L$ 
 with the group of all automorphims of $G$ that preserve $g$.   
\item"(ii)" 
Let $W=\{a\in G \mid R_a\text{ preserves }g\}$ and let 
$Z(G)$ be the {\rm(}finite{\rm)} center of~$G$. 
Then the identity component $\Isom_0(G,g)$ of the full isometry group 
is given by $(G_L\times W_R)/Z(G)$, where $Z(G)$ is embedded by 
$a\mapsto L_a\circ R_a\inv$. 
\endroster 
\endproclaim 
 
Note that $a\in W$ if and only if the associated inner automorphism $I_a$ 
of $G$ is an isometry of $(G,g)$. 
 
\definition{5.1.5 Notation} 
Let $P$ be the subgroup  of $G$ given by those matrices of the form 
$$\left(\matrix\Id_m&0&0\\0&\Id_m&0\\0&0&*\endmatrix\right). 
$$ 
Thus $P$ is isomorphic to $\SO(2r)$, respectively $\SU(r+1)$, in case 
$G=\somr$, respectively $G=\sumr$.   
 Let $Q=H\times P$. 
Recall that $H$ is the subgroup $\SO(m)\times\SO(m)$, 
respectively $\SU(m)\times\SU(m)$, of~$G$, embedded in the upper left 
corner, and that $T$ is the standard maximal torus in~$P$. 
\enddefinition 
 
\proclaim{5.1.6 Lemma} 
Let $\phi$ be defined as in~{\rm5.1.1}. 
Let $\Cal G$ denote the subgroup of $G_L\times G_R$ consisting of those 
elements which leave~$\phi$ invariant, and let $\Cal G_0$ denote the 
connected component of the identity in~$\Cal G$. Then: 
\roster 
\item"(i)" If $L_aR_b\in\Cal G$ then $a\cdot b\in Q$ and $a,b$ are elements 
of $\text{\rm S}(\O(m)\times\O(m)\times\O(2r))$, respectively of 
$\text{\rm S}(\U(m)\times\U(m)\times\U(r+1))$. In particular, $I_a$ and 
$I_b$ leave the subgroups $H$ and~$P$ of~$G$ invariant. 
\item"(ii)" 
$\Cal G_0=Q_L\times Q_R$ in case $G=\somr$, and\newline 
$\Cal G_0=\Inn_G(H)\times P_L\times T_R$ in case $G=\sumr$. 
\endroster 
\endproclaim 
 
\demo{Proof} 
(i) 
It is obvious from the definition of~$\phi$ that it assumes its maximum 
value precisely on $Q\subset G$. By determining first which elements 
of $G_L\times G_R$ preserve~$Q$, and then ruling out 
those of these candidates which obviously do not preserve~$\phi$ 
(taking into particular account that $c_1>c_2>0$) it is not hard 
to see that $a$ and $b$ must be of the asserted form. Moreover, 
$aQb=Q$ implies $a\cdot b\in Q$. 
 
(ii) 
If $G=\somr$ then (i) implies $\Cal G_0\subset Q_L\times Q_R$\,. 
Since $Q_L\times Q_R$ does leave~$\phi$ (as defined for this~$G$) 
invariant, the statement follows. 
 
Now let $G=\sumr$. In this case, (i) says that $\Cal G$, and hence 
$\Cal G_0$\,, is contained in 
$$Q_L\times Q_R\times\Inn_G(D)=H_L\times H_R\times P_L\times P_R 
  \times\Inn_G(D), 
$$ 
where $D$ is the subgroup of $G$ consisting of the matrices 
$\left(\smallmatrix \alpha\Id_m&&\\&\beta\Id_m&\\&&\gamma\Id_{r+1} 
\endsmallmatrix\right)\in G$ with $\alpha,\beta,\gamma\in e^{i\R}$. 
It is clear from the definition of~$\phi$ and the comments in~5.1.1 
that $\Inn_G(H)\times P_L\times T_R$ does preserve~$\phi$, hence is 
contained in~$\Cal G_0$\,. We have to show that it is in fact equal 
to~$\Cal G_0$\,. We claim that any $1$@-parameter family in $H_L\times 
P_R\times\Inn_G(D)$ which leaves~$\phi$ invariant is actually contained 
in~$T_R$\,; this will obviously imply the statement. Let $h(t),p(t),d(t)$ 
be $1$@-parameter families in $H$, $P$, and~$D$, respectively, 
with the property 
that $L_{h(t)}R_{p(t)}I_{d(t)}$ leaves~$\phi$ invariant. 
Write $h(t)=\left(\smallmatrix k_1(t)&&\\&k_2(t)&\\&&\Id\endsmallmatrix 
\right)$, where $k_1(t),k_2(t)$ are $1$@-parameter families in~$\SU(m)$, 
and write $d(t)$ as above with $1$@-parameter families 
$\alpha(t),\beta(t),\gamma(t)$ in $e^{i\R}$. Identify $p(t)$ with the 
corresponding $1$@-parameter family in $\SU(r+1)$. 
Note that the first two summands in the exponent of~$\phi$ are invariant 
under $L_{h(t)}R_{p(t)}I_{d(t)}$\,; consequently, this family must 
also leave the function $\psi:G\ni X\mapsto 
(\text{Re}\,d(A,C))^2+(\text{Re}\,d(E,F))^2\in\R$ invariant, 
where $X$ is written as in~5.1.1. 
 
We first show that the family $h(t)$ is trivial. Assume the contrary. 
Then $k_1(t)$ or $k_2(t)$ is nontrivial; without loss of generality 
we assume $k_1(t)$ is nontrivial. 
Choose~$t_0$ such that $k_1(t_0)$ has at least two different eigenvalues. 
Choose $A'\in\SU(m)$ such that $A'$ has less than~$m$ different eigenvalues, 
but $k_1(t_0)\inv A'$ has~$m$ pairwise different eigenvalues. Choose 
$0<c<1$. It is easy 
to see that there exist elements $X$ of~$G$, written as in~5.1.1, 
with $A=c k_1(t_0)\inv A'$ and $F=0$. It is also not hard to show that 
the condition that no column 
vector of~$C$ be orthogonal to any of the eigenvectors of~$A$ --- 
equivalently, $d(A,C)\ne0$ --- is a generic 
condition on this subset of~$G$, and that 
so is the condition $\text{Re}\,d(A,C) 
\ne0$. Choose $X\in G$ with these properties. 
Then $\psi(X)=(\text{Re}\,d(A,C))^2+0\ne0$, 
but $\psi(L_{h(t_0)}R_{p(t_0)}I_{d(t_0)}X)=(\text{Re}\,d(cA', 
k_1(t_0)Cp(t_0)\alpha(t_0)\gamma(t_0)\inv)^2+0=0$ since $A'$ has less than 
$m$ different eigenvalues (recall the comments in 5.1.1). 
This is a contradiction; so the family~$h(t)$ must indeed be trivial. 
 
Next we attack $p(t)$. 
Assume that the family $p(t)$ were not contained in the 
maximal torus~$T$ consisting of the diagonal elements in~$P$. 
Let $k$ be a diagonal matrix in $\SU(m)$ with $m$ different eigenvalues, 
and set $A:=ck$ for some $c$ with $0<c<1$. Choose~$t_0$ such that some 
column vector $p_i(t_0)$ of $p(t_0)$ has at least two nonzero entries. 
It is easy to see that there exist elements $X$ of~$G$, written as in~5.1.1, 
with the given~$A=ck$ and with $F=0$, such that some row vector
of~$C$ is orthogonal 
to $p_i(t_0)$. It is also not hard to show that the condition that $C$ 
have only nonvanishing entries --- equivalently: $d(A,C)\ne0$ --- 
is a generic condition on this subset of~$G$ (recall the choice of 
$p_i(t_0)$), and that so is the condition $\text{Re}\,d(A,C)\ne0$. 
Choose $X\in G$ with these properties. Then $\psi(X)=(\text{Re}\, 
d(A,C))^2\ne0$, but $\psi(L_{h(t_0)}R_{p(t_0)}I_{d(t_0)}X)=(\text{Re}\, 
d(A,Cp(t_0)\alpha(t_0)\gamma(t_0)\inv)^2=0$ since $Cp_i(t_0)$ has a zero 
entry and thus $d(A,Cp_i(t_0))=0$ (recall that $A$ is diagonal).
This is a contradiction, 
so the family $p(t)$ must be contained in~$T$. 
 
Finally, for every $X\in G$ with $F=0$ we have 
that $\psi(R_{p(t)}I_{d(t)}X)$ is equal to 
$(\text{Re}\,d(A,Cp(t)\alpha(t)\gamma(t)\inv)^2 
=(\text{Re}((\alpha(t)\gamma(t)\inv)^{m(r+1)}d(A,C))^2$ 
since $p(t)\in T$. 
The fact that for every given such~$X$ this expression must be constant 
in~$t$ implies $\alpha(t)=\gamma(t)$ 
for all~$t$. Analogously we obtain $\beta(t)=\gamma(t)$ for all~$t$, 
hence $d(t)=\alpha(t)\Id_{2m+r+1}$\,. But $d(t)\in\sumr$; hence the family 
$d(t)$ is trivial. 
\qed\enddemo 
 
\proclaim{5.1.7 Lemma} 
Assume that $j$ satisfies the genericity condition~{\rm5.1.2(a)}. 
Let $\Phi$ be any automorphism of~$Q$ which preserves the subgroups $H$ 
and~$P$ and which is an isometry of the metric~$\gl$\,.  Then $\Phi$ 
preserves~$T$, and $(\Phi_*\inv\restr\h\,,\Phi_*\restr\z)$ is a 
self-equivalence of~$j$ in the sense of Definition~{\rm3.1}. 
\endproclaim 
 
\demo{Proof} 
For $X\in\h$ and $U\in\p$ we have $\gl(X,U)=g_0(\la(X),U)$ since 
$\la\restr\p=0$ and since $\h$ is $g_0$@-orthogonal to~$\p$. 
{}From the facts that 
$\gl(X,U)=\gl(\Phi_*(X),\Phi_*(U))$ and that 
$\Phi_*$ preserves $\h$ and~$\p$, 
it follows that 
$$g_0(\la(X),U)=g_0(\la(\Phi_*X),\Phi_*U). 
$$ 
Using the genericity condition 5.1.2(a), we may conclude that 
$\Phi_*(\p\cap\z^\perp)\subset \p\cap\z^\perp$ and thus $\Phi_*(\z)=\z$. 
Letting $U\in\z$ in the displayed equation above, we then have 
$$g_0(X,j(U))=g_0\bigl(\Phi_*X,j(\Phi_*U)\bigr) 
  =g_0\bigl(X,\Phi_*\inv(j(\Phi_*U))\bigr), 
$$ 
and the lemma follows. 
\qed\enddemo 
 
\proclaim{5.1.8 Theorem} 
We use the notation of Example~{\rm3.7}. Let $\phi$ be the function 
given in~{\rm5.1.1} and let $\psi=\tau^*\phi$.  Assume that $j_1$ is not 
equivalent to~$j_2$ and that $j$ satisfies the genericity 
condition~{\rm5.1.2(a)}.  Then $\psi$ is not congruent to~$\phi$ under an 
isometry of $(G,\gl)$.   
\endproclaim 
 
Thus by Example 3.7, $\phi$ and $\psi$ are non-congruent isospectral 
potentials on $(G,\gl)$. 
 
\demo{Proof} 
Suppose $F$ is an isometry of $(G,\gl)$ such that $F^*\psi=\phi$. 
By Lemma~5.1.4, $F\in G_L\sdp\,\Aut(G)$. 
Recall that $\tau=I_c$\,, where 
$$c=\left(\smallmatrix 0&\Id&0\\ \Id&0&0\\0&0&\Id\endsmallmatrix\right). 
$$ 
Thus $I_c^*\psi=\phi$, so $\phi$ is invariant under $I_{c\inv}\circ F$. 
Let $\sigma$ be the outer automorphism of~$G$ given (in case 
$G=\SO(2m+2r)$) by conjugation by a diagonal matrix of determinant $-1$ 
all 
of whose diagonal entries lie in $\{\pm 1\}$, or (in case $G=\sumr$) by 
the map $A\mapsto\bar A$.  Observe that $\sigma$ leaves $\phi$ invariant. 
Since for both choices of~$G$, $\Inn(G)$ has index two in $\Aut(G)$, 
$I_c\inv\circ F$ may be written either as $\sigma\circ I_a\circ L_b$ 
or else as $I_a\circ L_b$ for some $a,b\in G$.  In either case 
$(I_a\circ L_b)^*\phi=\phi$, so Lemma~5.1.6(i) implies that 
$I_a$ and $I_b$ preserve $H$ and~$P$. 
Consequently, the automorphism 
$\Phi:=F\circ L_b\inv$ (which is either $I_{ca}$ or $I_c\circ\sigma\circ 
I_a$) is an isometry of $(G,\gl)$ which preserves both $H$ and $P$ and 
interchanges the two factors of $H$. 
By Lemma~5.1.7, we see that $(\Phi_*\inv\restr\h\,,\Phi_*\restr\z)$ is a 
self-equivalence of $j$ which intertwines $j_1$ and $j_2$\,, contradicting 
the fact that $j_1$ and $j_2$ are not equivalent. 
\qed\enddemo 
 
We now prepare for the proof that the conformally equivalent metrics in 
Example~3.7 are not isometric for our given choice of~$\phi$ and with the 
genericity hypotheses~5.1.2. 
 
\proclaim{5.1.9 Proposition} 
Suppose that $(M,g)$ is a compact Riemannian manifold admitting a 
free action 
by isometries by a torus~$T$ and that the orbits all have the same volume 
{\rm(}i.e., they are totally geodesic{\rm)}.  Let $\phi$ and $\psi$ be 
positive smooth $T$@-invariant functions on $M$.  Let $S$ be a non-trivial 
subtorus of~$T$ and let $F:(M,\phi g)\to (M,\psi g)$ be an isometry that 
carries $S$@-orbits to $S$@-orbits.  Then $F^*\psi=\phi$ and 
$F\in\Isom(M,g)$.  In particular {\rm(}letting $\psi=\phi${\rm)}, 
the normalizer of~$S$ in $\Isom(M,\phi g)$ is contained in 
$\{F\in\Isom(M,g)\mid F^*\phi=\phi\}$. 
\endproclaim 
 
\demo{Proof} 
The second statement is immediate from the first.  To prove that 
$F^*\psi=\phi$, first note that all $S$@-orbits have the same volume 
relative to the metric~$g$; in fact, they are totally geodesic.  Moreover, 
since $\phi$ and $\psi$ are $S$@-invariant, they are constant 
on each orbit. 
Since the isometry~$F$ preserves volumes, it must intertwine the values of 
$\phi$ and~$\psi$ on each orbit.  Hence $\psi\circ F=\phi$. Consequently, 
$F$ is an isometry of $(M,g)$. 
\qed\enddemo 
 
\proclaim{5.1.10 Lemma} 
Let {\rm5.1.2(a,b)} be satisfied. 
Denote by $C$  the connected component of the identity in $\{F\in 
\Isom(G,\gl)\mid F^*\phi=\phi\}$.  Then $C=Q_L\times T_R$ in case 
$G=\somr$, and $C=P_L\times T_R$ in case $G=\sumr$, 
where $Q$ and $P$ are as in Notation~{\rm5.1.5}. 
\endproclaim 
 
\demo{Proof} 
If $G=\somr$ then by Lemmas 5.1.4 and~5.1.6(ii) and the fact that $\phi$ 
and $\gl$ are invariant under $Q_L\times T_R$, 
we have $Q_L\times T_R\subset C\subset Q_L\times 
Q_R$\,. If $G=\sumr$ then the same lemmas and the fact that $\phi$ is 
invariant under $P_L\times T_R$ imply $P_L\times T_R\subset C 
\subset \Inn_G(H)\times P_L\times T_R$. 
By Lemma~5.1.7 and the fact that the group of automorphisms of the 
torus~$T$ is discrete, the only one-parameter subgroups of $\Inn_G(Q)$ 
which lie in 
$\Isom(\gl)$ (consequently also: the only such subgroups of $\Inn_G(H)$) 
are those which restrict to the identity on 
$\z$ and also on $j(\z)$. Since the centralizer of $j(\z)$ in~$\h$ 
is trivial by condition 5.1.2(b), and since the centralizer of~$\z$ 
in~$\p$ is~$\z$, these one-parameter subgroups must 
lie in $\Inn_G(T)$. The corollary follows. 
\qed\enddemo 
 
\proclaim{5.1.11 Lemma} 
Let $U$ be a compact Lie group and $T$ a torus in~$U$ of dimension greater 
than one.  Suppose that every subtorus of~$T$ of positive dimension has 
the 
same centralizer in~$U$.  Then $T$ lies in the center of~$U$ {\rm(}and 
thus its centralizer is~$U${\rm)}. 
\endproclaim 
 
\demo{Proof} 
Let $\uu$ denote the Lie algebra of $U$ and $\t\subset\uu$ the Lie algebra 
of~$T$.  Consider the simultaneous root space decomposition of the 
complexified space $\uu^\C$ with respect to the operators $\ad(X)$, 
$X\in \t$. 
If $Y=Y_1+iY_2$ (with $Y_1,Y_2\in\uu$) is a root vector, say 
$\ad(X)(Y)=i\alpha(X)Y$ for all $X\in\uu$, where $\alpha\in\t^*$, 
then either $\alpha\equiv 0$ or else $\ker(\alpha)$ is a subspace $\s$ 
of~$\t$ 
 of codimension one.  In the latter case the centralizer of the 
associated 
subtorus $S$ of~$T$ contains~$Y$ and thus 
properly contains the centralizer of~$T$, contradicting 
our hypothesis.  Thus $\alpha\equiv 0$ for every root $\alpha$. 
Hence $\t$ is 
central in~$\uu$ and $T$ is central in~$U$. 
\qed\enddemo 
 
\proclaim{5.1.12 Proposition} 
If the genericity conditions {\rm5.1.2(a,b)} are satisfied, then 
the identity component of the center of $\Isom_0(G,\phi\gl)$ is equal 
to~$T_R$\,. 
\endproclaim 
 
\demo{Proof} 
By Proposition~5.1.9 and Lemma~5.1.10, the centralizer in 
$\Isom_0(G,\phi\gl)$ 
of {\it any\/} nontrivial subtorus of~$T_R$ is given in case 
$G=\somr$ by $Q_L\times T_R$\, and in case 
$G=\sumr$ by $P_L\times T_R$\,. 
Thus by Lemma~5.1.11, $\Isom_0(G,\gl)$ equals $Q_L\times T_R$\,, 
respectively $P_L\times T_R$\,. 
Since in either case, $T_R$ is the identity component of the 
center of this group, the proposition follows. 
\qed\enddemo 
 
Of course, the proposition is also valid with $\phi$ replaced by $\psi$. 
 
\proclaim{5.1.13 Theorem} 
Let $\phi$ be the function given in~{\rm5.1.1} 
and let $\psi=\tau^*\phi$. Assume that $j_1$ is not equivalent to~$j_2$ 
and that $j$ satisfies the genericity conditions~{\rm5.1.2(a,b)}. 
Then the metrics $\phi\gl$ and $\psi\gl$ on $G$ are not isometric. 
\endproclaim 
 
\demo{Proof} 
Suppose that $F:(G,\phi\gl)\to (G,\psi\gl)$ were an isometry. 
By Proposition~5.1.12, the resulting isomorphism from $\Isom(G,\phi\gl)$ 
to 
$\Isom(G,\psi\gl)$ given by $\mu\mapsto F\circ\mu\circ F\inv$ carries 
$T_R$ 
to~$T_R$\,.  Hence $F$ carries $T_R$@-orbits in $G$ to 
$T_R$@-orbits. 
By Proposition~5.1.9, $F$ is an isometry of $\gl$ and satisfies 
$F^*\psi=\phi$, contradicting Theorem 5.1.8. 
\qed\enddemo 
 
\subheading{5.2 Nonisometry proof for the examples from 4.4/4.6}

\noindent
We prove here the nonisometry of the conformally equivalent metrics on
both the balls and the spheres in Examples~4.4/4.6.  As an aside, we
observe that the nonisometry of the metrics on the balls would follow
immediately from that on the boundary spheres.  However, we will not use
this observation below but instead prove the nonisometry of the metrics
on the balls directly, so that we will be able later (in~5.4) to adapt
the proof to address the nonisometry of the metrics in Remark~4.11,
where the boundary spheres may in fact be isometric.

In the context of Example~4.4 (resp\.~4.6), suppose that the linear
maps $j,j':\z=\R^2\to\so(m)$ (resp\. $c,c':\R^2\to\Sym_0(\R^3)$)
satisfy the conditions from Remark~4.5 (resp\. those specified in~4.6);
let $\nu,\nup$ be the associated $\R^2$@-valued $1$@-forms on~$\R^m$
(where $m=3$ in the case of~4.6).
Furthermore, assume that the positive smooth radial function~$\psi$
on~$\R^m$
satisfies the condition from~4.5 (resp\.~4.6); writing
$\psi(p)=\rho(|p|)$,
this means that $\rho\restr{[0,1]}$ attains its maximum precisely in~$0$.

For $M=B\subset\R^m\oplus\R^m\oplus\C^2$,
or $M=S\subset\R^m\oplus\R^m\oplus\C^2$, we are to show that the metrics
$\phi_1g_\la$ and $\phi_2g_\la$ on~$M$ are nonisometric,
where $\la$ is the $\R^2$@-valued $1$@-form associated with
$\nu\oplus\nup$ on~$\R^m\oplus\R^m$ as in~4.1(ii), and where the $\phi_i$
($i=1,2$) are defined as in Theorem~4.3.  By the proof of Theorem~4.3,
the metric $\phi_2g_\la$ is isometric to
$\phi_1g_{\lap}$, where $\lap$ is the $\R^2$@-valued
$1$@-form associated with
$\nup\oplus\nu$.  Following the strategy described in
Proposition~2.2 and Remark~2.3, it thus suffices to show
that the $1$@-forms $\la$ and $\lap$ satisfy the conditions \thetag{N}
and~\thetag{G} in Proposition~2.2 with respect to the metric~$\phi g_0$
on~$M$, where $\phi:=\phi_1$\,. Recall that $\phi_1(p,q,u)=\psi(p)
=\rho(|p|)\in\R_+$\,.

We canonically identify the set $\hat M/T$ (see Notation~1.2) with
$$B^{2m+2}\cap Q=\{(p,q,a,b)\in\R^m\oplus\R^m\oplus\R\oplus\R\mid
a,b>0, |p|^2+|q|^2+a^2+b^2<1\}
$$
in case $M=B$, resp\. with
$$S^{2m+1}\cap Q=\{(p,q,a,b)\in\R^m\oplus\R^m\oplus\R\oplus\R\mid
a,b>0, |p|^2+|q|^2+a^2+b^2=1\}
$$
in case $M=S$, where $Q$ denotes the open ``quadrant'' of~$\R^{2m+2}$
defined by $a,b>0$.
The metric $\phi g_0^T$ on~$\hat M/T$ is in each case given as
$\bar\phi g_0$\,, where $\bar\phi(p,q,a,b):=\rho(|p|)$, and where~$g_0$
now denotes the standard metric on~$\R^{2m+2}$ as well.

Let $\bar F\in\Autphbm$.

If $M=B$ then $\bar F$, being an isometry of $\bar\phi g_0$,
extends to a continuous self-map, again denoted~$\bar F$,
of the closure of $B^{2m+2}\cap Q$.
Note that $\bar F$ must preserve the $2m$@-dimensional ``edge''
$B^{2m}\times\{0\}\subset\R^{2m}\oplus\R^2$ because this is the image
under the projection $M\to M/T$ of those points whose $T$@-orbit is
just a point. Hence $\bar F$ must also preserve the opposite ``face''
$S^{2m+1}\cap Q$. So $\bar F$ preserves $S^{2m+1}\cap Q$ regardless of
whether
$M=B$ or $M=S$.

Note that the $T$@-orbit of a point $(p,q,v,w)\in M$ with $|v|=a$ and
$|w|=b$, endowed with the metric induced by~$\phi g_0$\,, is a rectangular
torus with side lengths $\sqrt{\rho(|p|)}\cdot 2\pi a$ and
$\sqrt{\rho(|p|)}\cdot 2\pi b$. Considering the ratio of the
side lengths, we see that $\bar F$ must either preserve or switch
the functions $(p,q,a,b)\mapsto a/b$ and $(p,q,a,b)\mapsto b/a$
on~$\hat M/T$ and thus on $S^{2m+1}\cap Q$ (by continuity in case $M=B$).

The length of the diagonal of the rectangular torus which is the
$T$@-orbit of a point $(p,q,a,b)\in\hat M/T$ is equal to
$\ell(p,q,a,b)$, where $\ell:\R^{2m+2}\ni(p,q,a,b)\mapsto
\sqrt{\rho(|p|)}\cdot 2\pi\sqrt{a^2+b^2}\in\R$.
Consequently $\bar F$ must preserve~$\ell$, too.
By our assumption on~$\psi$ (hence on~$\rho$), the restriction of~$\ell$
to $S^{2m+1}\cap Q$ attains its maximum precisely in the points
$(0,0,a,b)$; hence the set of these points, too, is invariant
under~$\bar F$.
Summarizing, for each vector $(a,b)\in S^1$ with $a,b>0$, the map~$\bar F$
either preserves or switches the two $2m$@-dimensional hemispheres
$S_{(a,b)}$ and $S_{(b,a)}$ where $S_{(a,b)}:=S^{2m+1}\cap
\{(p,q,ta,tb)\mid p,q\in\R^m,\,t\in\R_+\}$, and accordingly $\bar F$
fixes,
resp\. switches, their poles $(0,0,a,b)$ and
$(0,0,b,a)$.

Recall that $\bar F$ is, at the same time, an isometry
of~$\bar\phi g_0$\,.
In particular, it is a conformal self-map of the open domain
$B^{2m+2}\cap Q\subset\R^{2m+2}$ (in case $M=B$), resp\. of the open
domain
$S^{2m+1}\cap Q\subset S^{2m+1}$ (in case $M=S$).
By Liouville's theorem (see, e.g., \cite{Sp}, pp.~302--313),
and since the dimensions of these domains
are certainly greater than two, it follows that $\bar F$ is the
restriction
to~$\hat M/T$ of a M\"obius transformation~$\Phi$ of $\overline{\R^{
2m+2}}:=\R^{2m+2}\cup\{\infty\}$; that is, a composition of reflections
in hyperspheres and hyperplanes. (To see this for $S^{2m+1}\cap Q$,
one has to conjugate~$\bar F$ by the stereographical projection, which
is itself the restriction of a M\"obius transformation, and to use the
fact that each M\"obius transformation of $\overline{\R^{2m+1}}$ extends
to one of $\overline{\R^{2m+2}}$.)

{}From our previous arguments we now conclude that~$\Phi$, combined with
the reflection in the $(a=b)$@-hyperplane if necessary, preserves
each of the $2m$@-dimensional hemispheres $S_{(a,b)}$
considered above and, moreover, fixes their poles $(0,0,a,b)$.
By well-known properties of M\"obius transformations,
this implies that on each of these hemispheres,
$\Phi$ is uniquely determined by its restriction to the boundary
$(2m-1)$@-sphere, on which it acts as an orthogonal transformation. 
Since all the hemispheres $S_{(a,b)}$ have the same boundary, namely,
$S^{2m-1}\times\{0\}\subset\R^{2m}\oplus\R^2$, it follows that 
$\Phi\restr{S^{2m+1}}$
is of the form $(U,\sigma)$ with $U\in\O(2m)$ and $\sigma\in
\bigl\{\Id,\left(\smallmatrix 0&1\\1&0\endsmallmatrix\right)\bigr\}
\subset\O(2)$.

Recall once more that $\bar F$ --- and hence~$\Phi$ --- has to
preserve the function $\ell:(p,q,a,b)\mapsto\sqrt{\rho(|p|)}\cdot
2\pi\sqrt{a^2+b^2}$ on $S^{2m+1}\cap Q$. Since our assumption on~$\rho$
implies, in particular, that $\rho$ is nonconstant on $[0,1]$, it follows
that $U$ is of the form $(A,A')\in\O(m)\times\O(m)$.
Thus $\Phi\restr{S^{2m+1}}=(A,A',\sigma)$. If $M=B$ then $\Phi$~must
moreover
preserve~$B^{2m+2}$ because it extends~$\bar F$. Summarizing, we
have shown so far:
$$
\gathered
\text{Each $\bar F\in\Autphbm$ is the restriction to $\hat M/T\subset
\R^{2m+2}$ of}\\ \text{some map $(A,A',\sigma)\in\O(m)\times\O(m)\times
\O(2)$.}\endgathered\tag{1}
$$
Moreover, another look at the isometry classes of the $T$@-orbits
in $(M,\phi g_0)$ shows that each $\Psi\in\Cal D$ (see 2.1(ii))
will preserve the set $\{\pm Z_1\,,\pm Z_2\}\subset\z=\R^2$;
in particular, $\Psi\in\O(2)$.

\medskip

\subsubhead{Proof of condition~\thetag{N}}\endsubsubhead

Suppose condition~\thetag{N} were violated. Let $\Psi\in\Cal D$
and $\bar F\in\Autphbm$ be such that $\Om\subla=\Psi\circ\bar F^*
\Om\sublap$\,. Note that we have $d\om_0=0$ here and hence $\Om_0=0$;
thus $\Om\subla=d\bar\la$, $\Om\sublap=d\bar\lap$ with $\bar\la,\bar\lap$
as in~2.1(iv).
Let $\pr:\R^m\oplus\R^m\oplus\R^2\to\R^m\oplus\R^m$ denote the
canonical projection and also its restriction to $\hat M/T=
B^{2m+2}\cap Q$, resp\. $\hat M/T=S^{2m+1}\cap Q$.
Then
$$d\bar\la=\pr^*(d\nu\oplus d\nup),\quad d\bar\lap=
\pr^*(d\nup\oplus d\nu).
$$
For $X,Y\in T_p\R^m$ and $Z\in\z$, we have
$$\align
\<d\nu^{(\prime)}(X,Y),Z\>&=2\<j^{(\prime)}(Z)X,Y\>
\quad\text{(in~4.4), resp\.}\\
\<d\nu^{(\prime)}(X,Y),Z\>&=3\<c^{(\prime)}(Z)p\times X,Y\>,
\quad\text{(in~4.6)}
\endalign
$$
 (compare the proof of \cite{Sch3},~4.3, for the second equation).
By~\thetag{1}, $\bar F$ is the restriction to~$\hat M/T$ of some map
$(A,A',\sigma)\in\O(m)\times\O(m)\times\O(2)$. Moreover, the interior
of $\pr(\hat M/T)\in\R^{2m}$ is the whole open unit ball~$B^{2m}$,
regardless of whether $M=B$ or $M=S$.
Therefore, one elementarily derives from the equation
$d\bar\la=\Psi\circ\bar F^*d\bar\lap$ and the formulas above:
$$\aligned j'(\trans\Psi(Z))=Aj(Z)A\inv,\quad
&j(\trans\Psi(Z))=A'j'(Z)A^{\prime\,-1},\quad\text{resp\.}\\
c'(\det(A)\cdot\trans\Psi(Z))=Ac(Z)A\inv,\quad
&c(\det(A')\cdot\trans\Psi(Z))=A'c'(Z)A^{\prime\,-1}
\endaligned\tag{2}
$$
for all $Z\in\z$, where $\trans\Psi$ denotes the transpose of~$\Psi$
with respect to the standard inner product on $\z=\R^2$. Since
$\Psi\in\O(\z)$ (see above) and $\det(A),\det(A')\in\{\pm1\}$,
the equations in~\thetag{2} contradict the nonequivalence assumption
on $j,j'$ from~4.5, resp\. on $c,c'$ from~4.6.

\medskip

\subsubhead{Proof of condition~\thetag{G}}\endsubsubhead

Suppose to the contrary that $\lap$ did not have Property~\thetag{G};
i.e., there is a nontrivial $1$@-parameter family $\bar F_t\in\Autphbm$
such that $\bar F_t^*\Om\sublap\equiv\Om\sublap$\,.
Proceeding exactly as in the proof of condition~\thetag{N},
replacing~$\Psi$ by~$\Id$ and $\la$ by~$\lap$, we now obtain
$1$@-parameter
families $A_t\,,A'_t\in\SO(m)$ such that
$$\aligned j'(Z)\equiv A_t j'(Z)A_t\inv,\quad
&j(Z)\equiv A'_t j(Z)A_t^{\prime\,-1},\quad\text{resp\.}\\
c'(Z)\equiv A_t c'(Z)A_t\inv,\quad
&c(Z)\equiv A'_t c(Z)A_t^{\prime\,-1}
\endaligned\tag{3}
$$
for all $Z\in\z$ and $t\in\R$ (compare~\thetag{2}).
The families $A_t\,,A'_t$ are those occurring in
$\bar F_t=(A_t\,,A'_t\,,\Id)
\in\SO(m)\times\SO(m)\times\SO(2)$. Note that the $\SO(2)$@-component
of~$\bar F_t$ is trivial because $\bar F_t$ must preserve the
open quadrant~$Q$ of $\R^{2m+2}$ considered earlier.
Since the family $\bar F_t$ is assumed to be nontrivial, at least
one of the $1$@-parameter families $A_t\,,A'_t$ in $\SO(m)$
is nontrivial. The corresponding equations in~\thetag{3} thus
contradict the genericity assumptions from~4.5 on $j$ and~$j'$,
resp\. from~4.6 on $c$ and~$c'$.
\qed

\subheading{5.3 Nonisometry proof for Example 4.10}

\noindent
The proof will be similar to that of \cite{Sch3}, Proposition~4.3.
Let $c,c':\z=\R^2\to\Sym_0(\R^3)$ be the pair of linear maps
from Example~4.10, and let $\nu,\nup$ be the associated $\z$@-valued
$1$@-forms on~$\C^2$.
We are going to show that the $\z$@-valued $1$@-forms
$\la,\lap$ on $M=B^8\subset\C^2\oplus\C^2$ or
$M=S^7\subset\C^2\oplus\C^2$,
induced by $\nu\oplus\nup$, resp\. $\nup\oplus\nu$, satisfy conditions
\thetag{N} and~\thetag{G} from Proposition~2.2 with respect to the metric
$\phi g_0$ on~$M$, where $\phi:=\phi_1:\C^2\oplus\C^2\ni(p,q)
\mapsto\psi(p)
\in\R_+$\,.  Since the metrics $\phi_1g_{\lap}$ and $\phi_2g_\la$ on $M$
are isometric, as seen in the proof of Theorem 4.9, we may then conclude
from Proposition 2.2 that the conformally equivalent metrics
$\phi_1g_\la$ and $\phi_2g_\la$ are not isometric.

For $a,b\ge0$ let $M_{a,b}:=\{(p,q)\in\C^2\oplus\C^2\mid
|p|^2=a^2, |q|^2=b^2\}=S_a^3\times S_b^3$.
Then $\hat M$ is the union of those $M_{a,b}\subset M$ with $a,b>0$
and $a^2+b^2<1$ (in case $M=B^8$), resp\. $a^2+b^2=1$ (in case $M=S^7$).

Write $\psi(p)=\rho(|p|)$.
The $T$@-orbit of $(p,q)\in M_{a,b}$\,, endowed with the metric
induced by~$\phi g_0$\,, is isometric to the rectangular torus
with sides of length $\sqrt{\rho(a)}\cdot2\pi a$ and $\sqrt{\rho(a)}\cdot
2\pi b$.
Since $\rho$ is strictly increasing by assumption,
so is the function $a\mapsto
\sqrt{\rho(a)}\cdot a$; note that this function maps $(0,\infty)$
bijectively to itself.
Consequently, for every $(a,b)\in(0,\infty)^2$ there is exactly
one $(a',b')\in(0,\infty)^2$ such that
$$\sqrt{\rho(a')}\cdot a'=\sqrt{\rho(a)}\cdot b\quad\text{and}\quad
\sqrt{\rho(a')}\cdot b'=\sqrt{\rho(a)}\cdot a
$$
(the first condition determines~$a'$, the second then determines~$b'$).
The map $\sigma:(a,b)\mapsto(a',b')$ on $(0,\infty)^2$
is continuous and obviously has
the following property: Whenever $a>b$, then $a'<a$ (since
$a\mapsto\sqrt{\rho(a)}\cdot a$ is strictly increasing); hence
$\sqrt{\rho(a')}<\sqrt{\rho(a)}$ and thus $\|(a',b')\|^2>\|(a,b)\|^2$.
In particular, $\sigma(a,b)$ lies outside the closed unit disc for all
$(a,b)\in\Cal C:=\{(a,b)\in S^1\mid a>b>0\}$,
and the same is true for some open
neighborhood $\Cal A\subset\{(a,b)\mid a>b>0\}$ of this arc~$\Cal C$
of~$S^1$.
Choose $(a,b)\in\Cal C$ (in case $M=S^7$), resp\. $(a,b)\in\Cal U:=
\Cal A\cap
B_1(0)$ (in case $M=B^8$).
Then $M_{a,b}$ is precisely the set of points in~$(\hat M,\phi g_0)$
whose $T$@-orbit is isometric to the
rectangular torus with sides of length $\sqrt{\rho(a)}\cdot2\pi a$ and
$\sqrt{\rho(a)}\cdot2\pi b$.

Each $F\in\Autphm$ (see Definition~2.1)
must therefore preserve $M_{a,b}=:L$. Note that $(L/T,\phi g_0^T)$
is isometric to $(S^2_{a/2}\times S^2_{b/2}\,,\,\rho(a)^2
g_{\text{stdd}})$,
where~$g_{\text{stdd}}$ denotes the restriction of the standard metric
on $\R^3\oplus\R^3$.
Since $a>b$, every isometry of $(L/T,\phi g_0^T)$ ---
and thus the restriction to~$L/T$ of every $\bar F\in\Autphbm$ ---
must preserve the two factors. In particular, every
$F\in\Autphm$ must preserve
the factors of $M_{a,b}$\,; we conclude from this that every
$\Psi\in\Cal D$
is of the form $Z_1\mapsto\eps_1Z_1$, $Z_2\mapsto\eps_2Z_2$
with $\eps_1\,,\eps_2\in\{\pm 1\}$; in particular, $\Psi\in\O(\z)$.

\medskip

\subsubhead Proof of condition~\thetag{N}\endsubsubhead

Suppose condition~\thetag{N} were violated. Let $\Psi\in\Cal D$
and $\bar F\in\Autphbm$ be such that $\Om\subla=\Psi\circ
\bar F^*\Om\sublap$\,. Let $L=M_{a,b}\subset\hat M$ be chosen as above.
By the arguments above, $\bar F$ induces an isometry of $(L/T,\phi
g_0^T)$. Denote by $\Om\subla^L$\,, $\Om\sublap^L$ the $\R^2$@-valued
$2$@-forms induced by $\Om\subla$\,,
$\Om\sublap$ on~$L/T\subset\hat M/T$; then $\Om\subla^L=\Psi\circ\bar F^*
\Om\sublap^L$ because $\bar F$ preserves~$L/T$.

Observe that the first, resp\. second, component of $\Om_0^L$ is a
scalar multiple of the standard volume form of the first, resp\. second,
factor of~$L/T$ (compare~\cite{Sch3}, 4.1(iii)). Recall from~2.1(iv)
that $\Om\subla=\Om_0+d\bar\la$, $\Om\sublap=\Om_0+d\bar\la'$, and
similarly
for the forms induced on~$L/T$.
Note that each of $\Om_0^L$, $d\bar\la^L$ and $d\bar\la^{\prime\,L}$
is a sum of pullbacks of forms defined on one of the factors of~$L/T$.
Since each $2$@-form on~$S^2$ is uniquely decomposable into an
exact component and a multiple of the volume form, and since $\bar F$
preserves the two factors of~$L/T$, we therefore obtain
$d\bar\la^L=\Psi\circ\bar F^*d\bar\la^{\prime\,L}$.

The restriction of~$\bar F$ to~$L/T$ is given as $(A,A')$
with $A,A'\in\O(3)$.
Evaluating the previous equation on $2$@-vectors tangent to the first
factor of~$L/T$, we get
$d\bar\nu=\Psi\circ A^*d\bar\nu'$,
where
$$\<d\bar\nu^{(\prime)}(X,Y),Z\>=3a^3\<c^{(\prime)}(Z)x\times X,Y\>
$$
for all $X,Y\in T_xS_{a/2}^2$ and $Z\in\z$;
compare the proof of \cite{Sch3},~4.3.
As in~\cite{Sch3} one now derives
$$c'(\Phi(Z))=Ac(Z)A\inv\tag{4}
$$
for all $Z\in\z$, with $\Phi:=\det A\cdot\Psi\in\O(\z)$.
But this contradicts the nonequivalence of~$c$ and~$c'$
(see condition~2.)~in~4.6).

\medskip

\subsubhead Proof of condition~\thetag{G}\endsubsubhead

Suppose to the contrary that $\lap$ did not have Property~(G);
i.e., there is a nontrivial $1$@-parameter family $\bar F_t\in\Autphbm$
such that $\bar F_t^*\Om\sublap\equiv\Om\sublap$\,.
Note that $\bar F_t$\,, being an isometry of~$\hat M/T$, is determined
by its restriction to any nonempty open subset. The submanifolds
$L/T=M_{a,b}/T$ do foliate an open subset of $\hat M/T$ if we vary the
tuple
$(a,b)$ in the set from which we had chosen it.
It follows that $\bar F_t$ induces a nontrivial $1$@-parameter family
of isometries on at least one of these $M_{a,b}/T$;
again, we denote the corresponding~$M_{a,b}$ by~$L$.

Proceeding exactly as in the proof of condition~\thetag{N},
replacing $\Psi$ by~$\Id$ and $\la$ by~$\lap$,
we now obtain $1$@-parameter families $A_t\,,A'_t$ in~$\SO(3)$,
such that $d\bar\la^{\prime\,L}\equiv(A_t\,,A_t')^*d\bar\la^{\prime\,L}$
on~$L/T$. Evaluating this on the first, resp\. second, factor of~$L/T$,
we see that $d\bar\nu'$ is invariant under the~$A_t$\,,
resp\. that $d\bar\nu$ is invariant under the~$A'_t$\,.
We derive that
$$c'(Z)\equiv A_tc'(Z)A_t\inv\quad\text{and}\quad
  c(Z)\equiv A'_tc(Z)A_t^{\prime\,-1}
$$
for all $Z\in\z$ (compare~\thetag{4}).
However, since at least one of the families $A_t\,,A'_t$ is nontrivial,
this contradicts the genericity assumption on~$c$ and~$c'$
(see condition~3.)~in~4.6).
\qed
\subheading{5.4 Nonisometry proof for the modified metrics from
  Remark~4.11}

\subsubhead The case of Example 4.4/4.6\endsubsubhead

We have to suitably adapt the proof given in 5.2. Note that $\la,\lap$
did not play any role in the derivation of statement~\thetag{1} of~5.2.
In order to prove the conditions~\thetag{N} and~\thetag{G}
for~$f\la,f\lap$,
we proceed as follows:

Choose any tuple $(r,s,t)$ such that $\alpha(r,s,t)\ne0$; we can assume
$r,s,t\ne0$. Let $L:=\{(p,q,u)\in\hat M\mid|p|=r,|q|=s,|u|=t\}$.
Then $f$ has the constant value $C:=\alpha(r,s,t)\ne0$ on~$L$.
For all differential forms occurring in the proof of condition~\thetag{N}
in~5.2, we now restrict attention to the forms they induce on~$L$,
resp\. on~$L/T$. These induced forms are, up to multiplication by
the constant~$C$, the same as those induced by $f\la,f\lap$ etc.
The only problem in the course of the arguments from~5.2 is the fact
that $\pr(L/T)=S_r^{m-1}\times S_s^{m-1}$ has empty interior in~$\R^{2m}$.
But the equations in~\thetag{2} can be derived nevertheless, even when
using the equation $d\nu\oplus d\nup=\Psi\circ(A,A')^*(d\nup\oplus d\nu)$
only restricted to the tangent bundle of $S_r^{m-1}\times S_s^{m-1}$:

In the case of Example~4.4, this is obvious from the formula for $d\nu^{
(\prime)}$, because for all $X,Y\in\R^m$ there is certainly some
$p\in S_r^{m-1}$ such that $X,Y\in T_pS_r^{m-1}$ (and similarly for
$S_s^{m-1}$).
In the case of Example~4.6, we first obtain only
$$\bigl\langle\bigl(c(Z)-A\inv c'(\det(A)\cdot\trans\Psi(Z))A\bigr)p
  \times X,Y\bigr\rangle=0
$$
for all $Z\in\z$ and all $X,Y\in T_pS_r^2$ (and similarly for
$c,c'$ reversed, with $A'$ instead of~$A$).
But this implies $\bigl(c(Z)-A\inv c'(\det(A)\cdot\trans\Psi(Z))A\bigr)p
\perp p$ for all $p\in S_r^2$\,. Since on the left hand side a symmetric
map
is applied to~$p$, this map must be zero. Thus the equations
in~\thetag{2}
follow nevertheless.

For both examples, the proof of condition~\thetag{G} now goes through
verbatim.

\subsubhead The case of Example~4.10\endsubsubhead

In case $M=B^8$, we make the assumption that the support of~$\alpha$
nontrivially intersects~$\Cal U$, where $\Cal U\subset B_1(0)\cap
\{(a,b)\mid a>b>0\}$ is the nonempty open set occurring in the proof
given in~5.3. In case $M=S^7$, since $\alpha$ does not identically vanish
on the relevant parameter set $\{(a,b)\in S^1\mid a,b\ge0\}$, and because
of $\alpha(a,b)=\alpha(b,a)$, we know that $\alpha$ does not identically
vanish on the arc $\Cal C=\{(a,b)\in S^1\mid a>b>0\}$ either.

Thus we can choose $(a,b)\in\Cal U$, resp\. $(a,b)\in\Cal C$,
such that $C:=\alpha(a,b)\ne0$. In~5.3 we saw that every $\bar F\in
\Autphbm$ must then preserve $L:=M_{a,b}=S_a^3\times S_b^3\subset\hat M$.
Moreover, $f\restr L\equiv C$.
In order to prove the conditions~\thetag{N} and~\thetag{G} for
$f\la,f\lap$,
we now work with $C\la^L,C\la^{\prime\,L}$
instead of $\la^L,\la^{\prime\,L}$ in~5.3, which does not affect any
of the arguments.
\qed

\bigskip

\heading\S6 Other types of examples\endheading

\noindent
The techniques used in Sections 3 and 4 apply to the construction of
isospectral conformally equivalent metrics and isospectral potentials in
various settings.  We first observe the common properties of the
constructions, particularly the constructions of Section 3 and of
Example~4.4.  For simplicity, we will consider here linear maps
$j$ from a vector space $\z$ into $\so(m)$ rather than into a more
general compact Lie algebra $\h$.  

\item{(i)} Linear maps $j:\z\to\so(m)$ where $\z$ is a
$k$@-dimensional vector space, viewed as the Lie algebra of a torus~$T$,
were used to construct Riemannian metrics $g_j$ on a fixed
manifold $M=M_{k,m}$\,; the manifold $M$ depended only on $k$ and $m$.   
(The  metrics were denoted by $g_{\la^j}$ rather than
$g_j$\,.  In Section~3, $M$ was a Lie group, while in Section~4,
$M$ was either a ball or sphere.)
The torus~$T$ acted by  isometries on each $(M,g_j)$.  

\item{(ii)}  In each setting, we knew from previous work that
$g_j$ and $g_{j'}$ were isospectral whenever the linear maps $j$ and $j'$
were
isospectral in the sense of Definition~3.1.  Moreover, the isospectrality
of the metrics $g_j$ and $g_{j'}$ had been proven using Theorem~1.3 with
$\phi=\psi=0$.  Subtori $W\subset T$ of codimension one correspond to
linear functionals $\mu\in\LL^*$ or, equivalently (given an inner product
on $\z$) to elements $Z$ of $\z$.  The diffeomorphisms $F_W$  of
Theorem 1.3  (referred to in Sections~3 and~4 as $F_\mu$) arose from the
linear maps $A_Z$  in Defintion~3.1.

\item{(iii)}  Given isospectral linear maps $j_1,j_2:\z\to\so(m)$, we
defined isospectral maps $j=(j_1,j_2)$  and
$j'=(j_2,j_1):\z\to\so(2m)$.  The associated metrics $g_j$ and $g_{j'}$
on $M=M_{k,2m}$ were isometric.  Thus we had diffeomorphisms $F_W$
satisfying (ii), arising from the isospectrality of $j$ and $j'$, and
separately an isometry
$\tau$ between
$g_j$ and
$g_{j'}$.  We then chose a conformal factor or potential $\phi$ invariant
under all the $F_W$ but not under $\tau$.  The functions $\phi$ and
$\tau^*\phi$ then gave isospectral potentials or conformal factors on
$(M,g_j)$.

\medskip

\noindent
We now briefly describe how this procedure can be used to obtain
isospectral potentials and conformal factors in other settings.

\subheading{6.1 Example}
Let $T$ be a torus with Lie algebra $\z$. 
Let $N=N_m$ denote either the standard unit ball or the standard unit
sphere in $\R^m$.  Given $j:\z\to\so(m)$, define a $\z$-valued $1$@-form
$\la=\la_j$ on $\R^m$ by $\la_p(X)=\langle j(p),X\rangle$, where
$\scp$ is the standard inner product on $\R^m$, and where
$T_p\R^m$ is identified with $\R^m$.  Then $\la$ pulls back, first by
the inclusion, to a $1$@-form on $N$ and then, by the canonical
projection,
to an admissible form on $N\times T$.  Let $g_j$ denote the
associated Riemannian metric on $N\times T$ as in Notation~1.4.  As
seen in \cite{GW2} and
\cite{GGSWW}, isospectral linear maps $j,j':\z\to\so(m)$ result in
isospectral metrics $g_j$ and $g_{j'}$\,; in the notation of~(ii) above,
the role of the diffeomorphisms $F_W$ is played by the maps
$A_Z\times\Id:N\times T\to N\times T$.  In the notation of (iii), the
role of $\tau$ is played by $\tau_0\times\Id$ where $\tau_0$ is the
diffeomorphism of the ball or sphere $N_{2m}$ obtained by restriction of
the map $(x,y)\mapsto(y,x)$ of $\R^{2m}$.    Choose a function $\psi$ on
$N$ which is invariant under the action of $\O(m)\times\O(m)$ but which is
not $\tau_0$-invariant.  Pull $\psi$ back by the projection to $N\times
T$ to obtain a function $\phi$ satisfying the conditions of (iii). 

\medskip

\subheading{6.2 Example}
Let $T$ be a torus with Lie algebra $\z$, and let  $Q$ be a domain with
boundary in the hyperbolic space $H^{m+1}$ invariant under the
isometric action of $\O(m)$.
With respect to the upper half space model of $H^{m+1}$ with
coordinates $(a,x)$, $a\in\R^+$, $x\in\R^m$, the action of $\alpha\in
\O(m)$ on $H^{m+1}$ is given by $\alpha(a,x)=(a,\alpha(x))$.
The article \cite{GSz} constructed from
each linear map 
$j:\z\to\so(m)$ a
Riemannian  metric $g_j$ on $Q\times T$.  After multiplying the linear map
$j$ by a small constant if necessary, the metric
$g_j$ was shown to have negative curvature.  Isospectral linear maps
resulted in isospectral metrics, with the role of the diffeomorphisms in
(ii) being played by the restriction to $Q\times T$ of the maps
$A_Z\times\Id$ on $H^{m+1}\times T$.  The procedure outlined in (iii)
results in isospectral potentials or isospectral conformally equivalent
metrics on $Q\times T$ where now $Q$ is an $\O(2m)$-invariant domain in
$H^{2m+1}$.  The role of $\tau$ in (iii) is played by the
restriction to $Q\times T$ of $\tau_0\times\Id$, where $\tau_0\in\O(2m)$
is given as in Example 6.1.  To obtain a function
$\phi$ as in (iii), begin with a function~$\psi$ on $H^{2m+1}$ invariant
under the action of
$\O(m)\times\O(m)$ but not  $\tau_0$@-invariant.  The function $\phi$ can
then be taken to be the pullback to $Q\times T$ of $\psi\restr Q$\,.
As in (iii),  $\phi$ and
$\phi\circ \tau$ give isospectral potentials on the negatively curved
manifold $(Q\times T,g_j)$.  By choosing $\phi$ to be positive, one
obtains isospectral conformally equivalent metrics on $Q\times T$.  By
moreover choosing
$\phi$ to be bounded by a sufficiently small constant, one can guarantee
that the conformally equivalent metrics are negatively curved.

If we choose $m=4\ell$ with $\ell\geq 2$, choose $\z$ to be 3-dimensional
and
choose $j_1\,,j_2: \z\to\so(m)$ as in Theorem~3.8(i) of \cite{GSz},
then the
metric $g_j$ on $Q\times T$, where $Q$ is a domain in $H^{2m+1}$ and
$j=(j_1\,,j_2)$, has constant Ricci curvature; indeed $(Q\times T,g_j)$
is a
domain in a harmonic manifold.  We thus obtain isospectral potentials  on
an Einstein manifold with boundary.

\enddocument